\pdfoutput=1
\documentclass[a4paper,english]{article}
\usepackage[T1]{fontenc}
\usepackage[latin9]{inputenc}
\usepackage{color}
\usepackage{babel}
\usepackage{float}
\usepackage{textcomp}
\usepackage{amsthm}
\usepackage{amsmath}
\usepackage{amssymb}
\usepackage{graphicx}
\usepackage[unicode=true,
 bookmarks=true,bookmarksnumbered=true,bookmarksopen=true,bookmarksopenlevel=4,
 breaklinks=false,pdfborder={0 0 1},backref=false,colorlinks=true,pdfpagemode=FullScreen]
 {hyperref}
\hypersetup{pdftitle={Latin Puzzles},
 pdfauthor={Miguel G. Palomo},
 pdfsubject={Custom Sudoku, Sudoku Ripeto, Combinatorics},
 pdfkeywords={Custom Sudoku, Sudoku Ripeto, Combinatorics},
 linkcolor=blue,citecolor=blue,filecolor=blue,urlcolor=blue}

\makeatletter

\pdfpageheight\paperheight
\pdfpagewidth\paperwidth

\providecommand{\tabularnewline}{\\}
\newcommand{\lyxdot}{.}

\floatstyle{boxed}
\newfloat{algorithm}{tbp}{lol}
\providecommand{\algorithmname}{Algorithm}
\floatname{algorithm}{\protect\algorithmname}

\numberwithin{equation}{section}
\numberwithin{table}{section}
\theoremstyle{plain}
\newtheorem{thm}{\protect\theoremname}[section]
\theoremstyle{definition}
\newtheorem{defn}[thm]{\protect\definitionname}
\theoremstyle{definition}
\newtheorem{example}[thm]{\protect\examplename}
\ifx\proof\undefined
\newenvironment{proof}[1][\protect\proofname]{\par
\normalfont\topsep6\p@\@plus6\p@\relax
\trivlist
\itemindent\parindent
\item[\hskip\labelsep\scshape #1]\ignorespaces
}{%
\endtrivlist\@endpefalse
}
\providecommand{\proofname}{Proof}
\fi

\date{}

\renewcommand\paragraph{\@startsection{paragraph}{4}{\z@}%
   {-3.25ex\@plus -1ex \@minus -.2ex}%
   {.75ex \@plus .2ex}%
   {\normalfont\normalsize\bfseries}}
\usepackage{caption}
\usepackage[margin=10pt,font=small,labelfont=bf,labelsep=period]{caption}

\AtBeginDocument{

}

\makeatother

\usepackage{listings}
\providecommand{\definitionname}{Definition}
\providecommand{\examplename}{Example}
\providecommand{\theoremname}{Theorem}

\begin{document}

\title{Latin Puzzles}

\author{Miguel G. Palomo}
\maketitle
\begin{abstract}
\noindent \addcontentsline{toc}{section}{Abstract}Based on a previous
generalization by the author of Latin squares to Latin boards, this
paper generalizes partial Latin squares and related objects like partial
Latin squares, completable partial Latin squares and Latin square
puzzles. The latter challenge players to complete partial Latin squares,
\textit{Sudoku} being the most popular variant nowadays.

The present generalization results in \emph{partial Latin boards},
\emph{completable partial Latin boards} and \emph{Latin puzzles}.
Provided examples of Latin puzzles illustrate how they differ from
puzzles based on Latin squares. The examples include\emph{ Sudoku
Ripeto} and \emph{Custom Sudoku}, two new \textit{Sudoku} variants.
This is followed by a discussion of methods to find Latin boards and
Latin puzzles amenable to being solved by human players, with an emphasis
on those based on constraint programming. The paper also includes
an analysis of objective and subjective ways to measure the difficulty
of Latin puzzles.

\textit{Sudoku}, \textit{Sudoku Ripeto}.
\end{abstract}

\section{Introduction\label{sec:Introduction}}

\begin{figure}
\begin{centering}
\includegraphics[height=4.5cm]{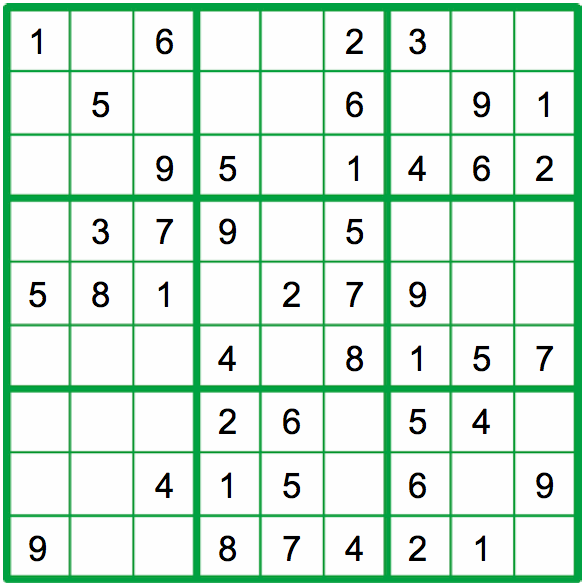}~~~~~~\includegraphics[height=4.5cm]{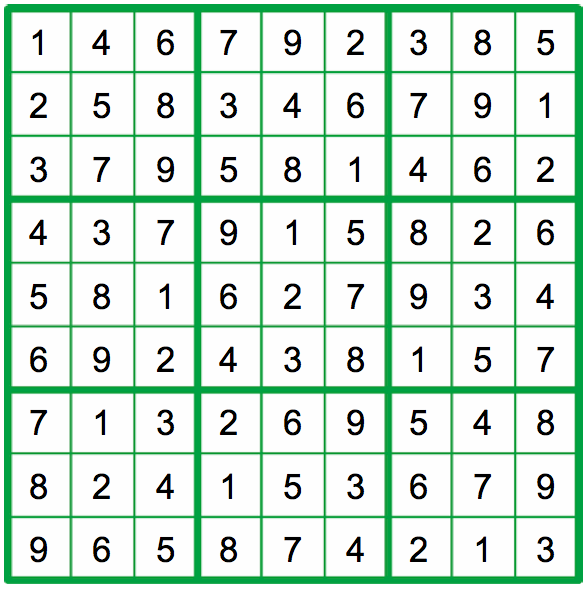}
\par\end{centering}

\protect\caption{\label{fig:Sudoku puzzle}Sudoku}
\end{figure}

\emph{Sudoku} puzzles challenge players to complete a square board
so that every row, column and $3\times3$ sub-square contains all
numbers from $1$ to $9$ (see an example in Fig. \ref{fig:Sudoku puzzle}).
The simplicity of the instructions coupled with the entailed combinatorial
properties have made \emph{Sudoku} both a popular puzzle and an object
of active mathematical research.

\noindent Back in the $19$\textsuperscript{th} century, B. Meyniel
created puzzles very similar to \textit{Sudoku} that were featured
in French newspapers \cite{Sudoku's French ancestors}\emph{. }The
puzzle was rendered popular in Japan in the late $1980$s by puzzle
editor Maki Kaji, who coined the word ``Sudoku'' (\emph{single number}
in Japanese). Later, former Hong Kong judge Wayne Gould wrote a computer
program to generate \emph{Sudoku} puzzles and had then published in
European and American newspapers. By the mid-2000s the puzzle was
popular worldwide. It was later pointed out by puzzle author and editor
Will Shortz that the creator of \emph{Sudoku} was American architect
Howard Garns ($1905$-$1989$), whose puzzles first appeared in 1979
in America under the name \emph{Number Place} \cite{Howard Garns}. 

Based on the generalization of Latin squares to Latin boards proposed
by the author in \cite{Latin Polytopes arxiv}, this paper proposes
a broad generalization of \emph{Sudoku} regarding board topology,
regions, symbols and inscriptions inside the board. This results in
new categories of puzzles in addition to the existing one, to which
\emph{Sudoku} and current variants belong.

\section{The Latin Square Puzzle\label{sec:The-Latin-Square-Puzzle}}

\subsection{Latin Squares}

\label{def:Latin-square}
\begin{figure}
\centering{}\includegraphics[height=4.5cm]{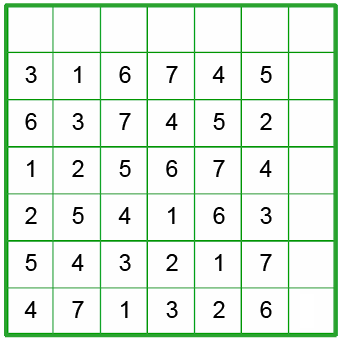}~~~~~~~\includegraphics[height=4.5cm]{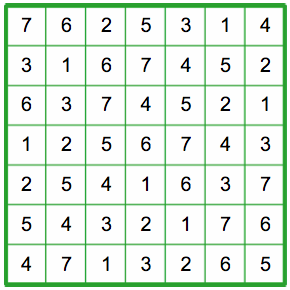}\protect\caption{\label{fig:LatinSquareOrder9}Latin square puzzle and solution}
\end{figure}

A Latin square of order \emph{n} is an $n\times n$
array of cells in which every row and every column holds all labels
belonging to an \emph{n}-set \cite{Handbook of Combinatorial Designs}.
Latin squares are so named because the $18$\textsuperscript{th}
century Swiss mathematician Leonhard Euler used Latin letters as labels
in his paper \emph{On Magic Squares} \cite{De quadratis magicis}.
Fig. \ref{fig:LatinSquareOrder9} (right) shows a Latin square of
order seven with the set of labels $\left\{ 1,2,3,4,5,6,7\right\} $.

Latin squares are redundant structures: after removing some labels
the result may still be completable to just the original Latin square.
If presented with the \emph{challenge}: ``Complete the board in Fig.
\ref{fig:LatinSquareOrder9} (left) so that every row and every column
contains all numbers from $1$ to $7$'', we can easily conclude
that the \emph{only} solution is the Latin square on its right\footnote{For some historical facts about Latin squares including ludic ones
see \cite{The history of Latin squares Andersen,An illustrated philatelic introduction,Some comments on Latin squares and on Graeco-Latin squares}.}. 

\noindent 

\subsection{\noindent Partial Objects Related to Latin squares\label{sub:Partial-Objects-Related-to-LS}}

We list here definitions of classic objects related to Latin squares
(see \cite{Latin squares and their applications}).\medskip{}

\begin{figure}
\begin{centering}
\includegraphics[height=4.5cm]{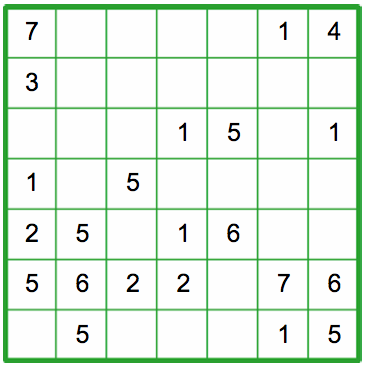}~~~~~~\includegraphics[height=4.5cm]{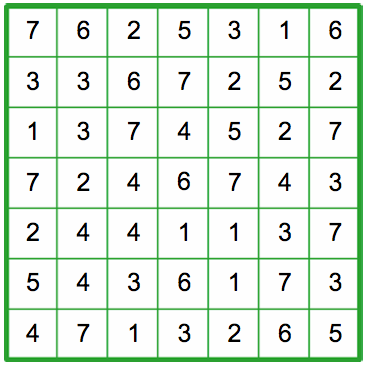}
\par\end{centering}

\protect\caption{\label{fig:partial magma and magma}Partial magma and magma}
\end{figure}

\noindent \textbf{Magmas}\label{sub: def magma}. A \emph{magma} of
order $n$ is a square grid of cells labeled with the elements of
an $n$-set of labels\footnote{\noindent There is a more algebraical definition: ``A \emph{magma}
is a set with an operation that sends any two elements in the set
to another element in the set''. We use then the term in this paper
to denote a magma's multiplication table without its borders.}. Fig. \ref{fig:partial magma and magma} (right) is an example with
labels $\left\{ 1,2,3,4,5,6,7\right\} $, some of them repeated in
rows and columns. \bigskip{}

\noindent \textbf{Partial Magmas}\label{Partial-magmas}. A \emph{partial
magma} is a magma with some labels missing. Fig. \ref{fig:partial magma and magma}
(left) shows an example with labels $\left\{ 1,2,3,4,5,6,7\right\} $.\bigskip{}

\begin{figure}
\begin{centering}
\includegraphics[height=4.5cm]{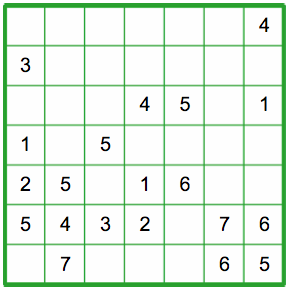}~~~~~~\includegraphics[height=4.5cm]{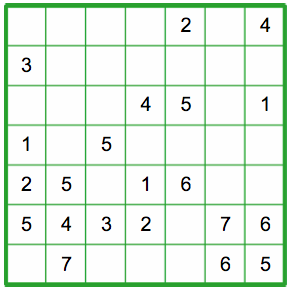}\medskip{}

\par\end{centering}

\begin{centering}
\includegraphics[height=4.5cm]{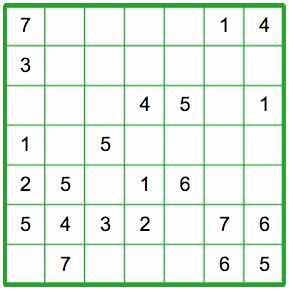}~~~~~~\includegraphics[height=4.5cm]{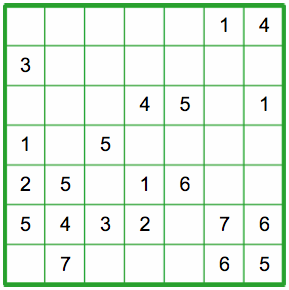}
\par\end{centering}

\protect\caption{\label{fig:PartialLatinSquares}Partial Latin squares}
\end{figure}

\noindent \textbf{Partial Latin Squares}\label{Partial-Latin-Squares Def}.
Fig. \ref{fig:PartialLatinSquares} shows four partial magmas with
labels from $\left\{ 1,2,3,4,5,6,7\right\} $ in which no label occurs
more than once in any row or column. Partial magmas fulfilling this
condition are called \emph{partial Latin squares} \cite{Latin squares and their applications}.\bigskip{}

\begin{figure}
\begin{centering}
\includegraphics[height=4.5cm]{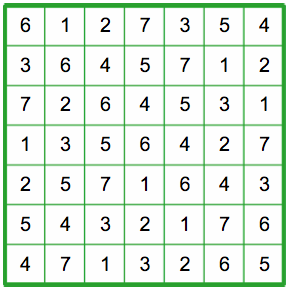}~~~~~~\includegraphics[height=4.5cm]{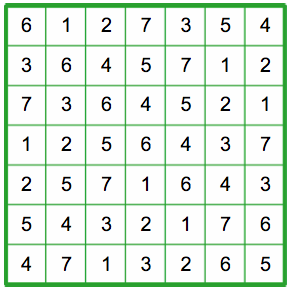}\medskip{}

\par\end{centering}

\begin{centering}
\includegraphics[height=4.5cm]{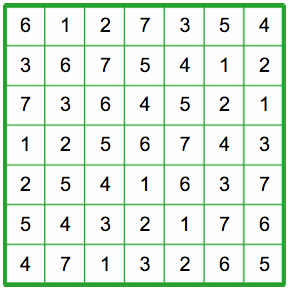}~~~~~~\includegraphics[height=4.5cm]{maxDifficulty_Less1j_solution_4of4_alsoSolOfBothSubAndCriticalSets}
\par\end{centering}

\protect\caption{\label{fig:solutions4partialLatinSquares}Latin squares}
\end{figure}

\noindent \textbf{Completable Partial Latin Squares}. We may wonder
if each of the squares in Fig. \ref{fig:PartialLatinSquares} can
be completable to either one or several Latin squares. It turns out
that the top-left one can only be completable to any of the four Latin
squares in Fig. \ref{fig:solutions4partialLatinSquares}, while the
top-right one can't be completed to any at all. This is easy to see
provided that the claim for the top-left square is true: the top-right
one only differs from the former in the extra $2$ on the first row,
but none of the four squares in Fig. \ref{fig:solutions4partialLatinSquares}
has $2$ in that position.

The two bottom squares in Fig. \ref{fig:PartialLatinSquares} can
be completable to \emph{just} the bottom-right square in Fig. \ref{fig:solutions4partialLatinSquares}.
We say that all partial Latin squares in Fig. \ref{fig:PartialLatinSquares},
except for the top-right one, are \emph{completable partial Latin
squares}.\bigskip{}

\noindent \textbf{Uniquely completable partial Latin squares}\label{Uniquely completable partial Latin squares}.
A partial Latin square that is completable to just one Latin square
\textendash like the two bottom ones in Fig. \ref{fig:PartialLatinSquares}\textendash{}
is a \emph{uniquely completable partial Latin square}.\bigskip{}

\noindent \textbf{Critical sets}. If the removal of any label in a
uniquely completable partial Latin square destroys its property of
uniqueness of completion \textendash the case of the bottom-right
square in Fig. \ref{fig:PartialLatinSquares}\textendash{} then the
square is a \emph{critical set}.\bigskip{}

\noindent \textbf{Latin square puzzles} \textbf{and} \textbf{critical
Latin square puzzles}. Uniquely completable partial Latin squares
and critical sets are called respectively \emph{Latin square puzzles\label{Latin square puzzle definition}}
and \emph{critical} \emph{Latin square puzzles} in this paper.

\section{Latin Puzzles}

We reproduce for convenience in this section some definitions found
in \cite{Latin Polytopes arxiv}.

\subsection{Boards}
\begin{defn}
\label{def:board}A \textit{board }is a pair $(P,C)$ where \textit{$P$}
is a set of \textit{\emph{geometric points}} and \textit{$C$} is
a set of subsets of $P$ whose union is\textit{\emph{ }}\textit{$P$.}
\end{defn}
\noindent We deduce from the definition that intersection among subsets
is allowed. The set $C$ is the \emph{constellation} of the board,
and every subset an \emph{asterism}.
\begin{defn}
\label{def:uniform board}A \emph{k-uniform board} is one whose asterisms
have all size $k$.
\end{defn}

\subsubsection{Points and Cells in Boards}

\noindent We will often represent a board with a grid of cells. To
conciliate this with the formal definition, that uses points, we adopt
the following convention: ``cell'' will mean either \emph{centroid
of the cell} or \emph{point,} depending on the context. For example
``write a label in the cell'' will mean \emph{label the centroid
of the cell}, whereas ``a board with $n$ cells'' will be \emph{a
board with $n$ points}.
\begin{example}
\noindent If we remove the numbers inside the cells in Fig. \ref{fig:LatinSquareOrder9}
(left), the resulting object is a board whose points are the (centroids
of) cells. Every row and column of (centroids of) cells is an asterism
of this board.
\end{example}

\subsection{Latin Boards}
\begin{defn}
\label{def:Latin board definition}Let $B=\left(P,C\right)$ be a
$k$-uniform board and \emph{L }a $k$-multiset of \emph{labels}.
A\emph{ Latin board} is a tuple $\left(B,L,F\right)$ where $F$ is
a bijective function $P\rightarrow L$.
\end{defn}
\noindent In simpler terms: in a Latin board every asterism holds
all labels in the multiset\footnote{This definition is for the so-called \emph{k-unifit Latin board}s.
See \cite{Latin Polytopes arxiv} for a more general definition of
Latin boards.}. We also say that a Latin board is the result of \emph{latinizing}
or \emph{labeling} a board with a multiset \cite{Latin Polytopes arxiv}.
According to this definition, and in contrast with the one for Latin
squares in Def. \ref{def:Latin-square}, Latin boards have no constraints
on topology, number of asterisms or uniqueness of labels.

\subsection{Partial Objects Related to Latin Boards}

We define here for the first time partial objects related to Latin
boards. We do this by generalizing partial objects related to Latin
squares (see Sect. \ref{sub:Partial-Objects-Related-to-LS}).

\subsubsection{Boards}

As per Def. \ref{def:board}, boards generalize the underlying $n\times n$
array of cells in Latin squares and their distribution in rows and
columns. Every row and column generalizes to an asterism; the set
of all rows and columns generalizes to the constellation of the board.

\subsubsection{Labeled Boards}
\begin{defn}
\noindent Let $B=\left(P,C\right)$ be a $k$-uniform board and \emph{L
}a $k$-multiset of \emph{labels}. \emph{A $k$-labeled board} is
a tuple $\left(B,L,F\right)$ where $F$ is a function $P\rightarrow L$.
\end{defn}
\noindent As the function needs not be bijective, the count of each
label in an asterism may not match its count in the multiset. Labeled
boards generalize magmas.

\subsubsection{Partial Labeled Boards\label{sub:def Partial-Labeled-Board}}
\begin{defn}
\noindent A\emph{ partial labeled board} is a labeled board with some
missing labels.
\end{defn}
\noindent Partial labeled boards generalize partial magmas.

\subsubsection{Partial Latin Boards}
\begin{defn}
\noindent \label{def:partial Latin board} A partial labeled board
is a\emph{ partial Latin board} if and only if the count of each label
in every asterism is not greater than its count in the multiset.
\end{defn}
\noindent Partial Latin boards generalize partial Latin squares.

\subsubsection{Completable Partial Latin Boards\label{sub:Completable-Partial-Latin-Board}}

\noindent As with partial Latin squares, a partial Latin board is
not guaranteed to be completable, hence the following definition.
\begin{defn}
A partial Latin board is a \textit{completable}\emph{ partial Latin
board} if and only if it can be completed to one or several Latin
boards.
\end{defn}
\noindent Completable partial Latin boards generalize completable
partial Latin squares.

\subsubsection{Latin puzzles}

\noindent In order to remove ambiguity when possible, we will follow
hereafter the\emph{ puzzle-Puzzle name convention}:\textit{ }\textit{\emph{``a
puzzle''}} (lower-case \emph{p}) is an instance of a particular ``Puzzle''
(upper-case \emph{p}), which is the general game \textendash ``\textit{\emph{a}}\textit{
Sudoku} puzzle'' vs. ``\textit{\emph{the}}\textit{ Sudoku} Puzzle''
for example.
\begin{defn}
\label{def:Latin-puzzle}A completable partial Latin board is a\textit{
Latin puzzle} if and only if it can be completed to exactly one Latin
board.
\end{defn}
\noindent We call this Latin board the puzzle's \emph{solution}. The
labels in the puzzle are called \emph{clues}. Latin puzzles generalize
uniquely completable partial Latin squares. 

Like Latin squares, Latin boards are redundant structures: they contain
more information than is strictly necessary to describe them, hence
the possibility of Latin puzzles. 

\noindent The board in Fig. \ref{fig:LatinSquareOrder9} (left) is
then a Latin puzzle, and also the two bottom ones in Fig. \ref{fig:PartialLatinSquares}.
The two top ones in the same figure are not Latin puzzles. Latin puzzles
generalize Latin square puzzles in the sense that their solutions
are Latin boards, themselves a generalization of Latin squares \cite{Latin Polytopes arxiv}.

\subsubsection{Latin Puzzles}
\begin{defn}
A \emph{Latin Puzzle} is the set of puzzles that results when we render
variable both the multiset of labels and the set of clues in a given
Latin puzzle.
\end{defn}
\noindent We will also say that a particular Latin puzzle is a \emph{member}
or an \emph{instance} of the corresponding Latin Puzzle. The three
Latin puzzles in Fig. \ref{fig:PartialLatinSquares} belong then to
the \emph{Latin Square Puzzle}. Latin Puzzles generalize the \emph{Latin
square Puzzle}.

\subsubsection{\noindent Inscripted Latin Puzzles}

\noindent Multisets of labels open the door to Puzzles with inscriptions.
\begin{defn}
\noindent An \emph{inscription} in a Latin Puzzle is a set of pairs
(label, cell) always present among the clues of every Puzzle's puzzle.
\end{defn}
\noindent We say that both the Puzzle and his puzzles are \emph{inscripted}.
Inscriptions make possible puzzles with generic words in any language
and alphabet, see Sect. \ref{sub:Inscriptions} and \cite{Latin Puzzles site}
for examples.

\subsubsection{The Latin Puzzle}
\begin{defn}
\emph{The Latin Puzzle} is the set of puzzles that results when we
render variable the board, the multiset of labels and the clues in
a given puzzle.
\end{defn}

\subsubsection{Critical Latin puzzles}
\begin{defn}
\label{def:critical-Latin-puzzle}A Latin puzzle is \textit{critical}
if the removal of any of its clues destroys the uniqueness of the
completion. An inscripted Latin puzzle is \textit{critical} if the
removal of any of its clues not in the inscription destroys the uniqueness
of the completion (see Sect. \ref{Critical Custom-Sudoku}).
\end{defn}

\subsubsection{Minimal Latin puzzles}
\begin{defn}
A Latin puzzle is \textit{minimal} if it has the smallest possible
set of clues among the puzzles in its Latin Puzzle.
\end{defn}

\subsubsection{Corresponding Terms}

Table \ref{tab:corresponding terms for LSs and LBs} summarizes the
correspondences between terms for Latin squares and Latin boards.

\begin{table}
\begin{raggedright}
\begin{tabular}{rl}
\noalign{\vskip0.1cm}
\textbf{\footnotesize{}Latin Square Term} & \textbf{\footnotesize{}Latin Board Term}\tabularnewline
\noalign{\vskip0.1cm}
{\footnotesize{}$n\times n$ array of cells} & {\footnotesize{}set of geometric points}\tabularnewline
\noalign{\vskip0.1cm}
{\footnotesize{} cells grouped in rows and columns} & {\footnotesize{}board}\tabularnewline
\noalign{\vskip0.1cm}
{\footnotesize{}row or column of cells} & {\footnotesize{}asterism}\tabularnewline
\noalign{\vskip0.1cm}
{\footnotesize{}set of all rows and columns} & {\footnotesize{}constellation}\tabularnewline
\noalign{\vskip0.1cm}
{\footnotesize{}magma} & {\footnotesize{}labeled board}\tabularnewline
\noalign{\vskip0.1cm}
{\footnotesize{}partial magma} & {\footnotesize{}partial labeled board}\tabularnewline
\noalign{\vskip0.1cm}
{\footnotesize{}partial Latin square} & {\footnotesize{}partial Latin board}\tabularnewline
\noalign{\vskip0.1cm}
{\footnotesize{}completable partial Latin square} & {\footnotesize{}completable partial Latin board}\tabularnewline
\noalign{\vskip0.1cm}
{\footnotesize{}uniquely completable partial Latin square} & {\footnotesize{}Latin puzzle}\tabularnewline
\noalign{\vskip0.1cm}
{\footnotesize{}critical set} & {\footnotesize{}critical Latin puzzle}\tabularnewline
\noalign{\vskip0.1cm}
\emph{\footnotesize{}Latin square Puzzle} & {\footnotesize{}Latin Puzzle}\tabularnewline
\noalign{\vskip0.1cm}
{\footnotesize{}Latin square} & {\footnotesize{}Latin board}\tabularnewline
\end{tabular}
\par\end{raggedright}

\protect\caption{\label{tab:corresponding terms for LSs and LBs}Correspondence between
classic Latin squares terms and new Latin boards ones}

\end{table}

\subsection{Types of Latin Puzzles\label{sub:Types-of-Latin-Puzzles}}

In the light of Defs. \ref{def:Latin board definition} and \ref{def:Latin-puzzle}
we may distinguish these main types of Latin Puzzles:
\begin{itemize}
\item \emph{LS} if their solutions are Latin squares, \emph{non-LS} otherwise
\item \emph{symmetric} if their solutions are Latin polytopes \cite{Latin Polytopes arxiv},
\emph{non-symmetric} otherwise
\item \emph{repeat} if the multiset has repeated labels, \emph{non-repeat}
otherwise
\item \emph{inscripted} if every puzzle in the Puzzle has an inscription
(see Def. \ref{sub:Inscriptions}), \emph{non-inscripted} otherwise
\end{itemize}
The \textit{Sudoku} Puzzle for example is then LS, symmetric, non-repeat,
non-inscripted.

\section{Latin Puzzles Examples\label{sec:Examples of Latin puzzles}}

The Latin puzzles and boards that follow have been created with the
methods described in Sect. \ref{sec:Finding-Latin-puzzles}. Each
puzzle is identified in bold typeface by its board name plus an adjective
indicating whether or not inscriptions, sets or multisets are used.
Except for the examples in Sect. \ref{sub:Latin-polyhedra}, they
are all guaranteed to be solvable in polynomial time in the number
of cells (see Sect. \ref{sub:Fair-Latin-puzzles}).

Each puzzle is included in a section that illustrates a specific relaxation
of the Latin square Puzzle (see Sect. \ref{sec:The-Latin-Square-Puzzle}).
Unless otherwise stated, each puzzle is rated \emph{medium} by the
method described in Sect. \ref{sub:Rating fair Latin puzzles}. Additional
Latin boards and puzzles can be found in \cite{Latin Polytopes arxiv,Latin Puzzles site}.

\subsection{Extra Asterisms}

By adding extra asterisms to the Latin square Puzzle we obtain \emph{Latin
square Puzzle variants}: those whose instances' solutions are Latin
squares.\\

\begin{figure}
\begin{centering}
\includegraphics[height=4.5cm]{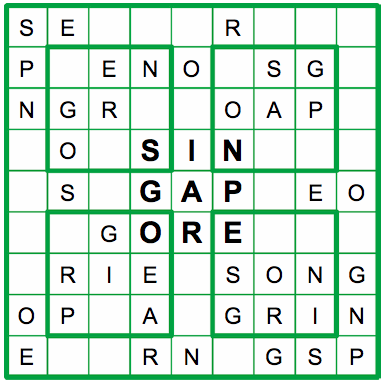}~~~~~~\includegraphics[height=4.5cm]{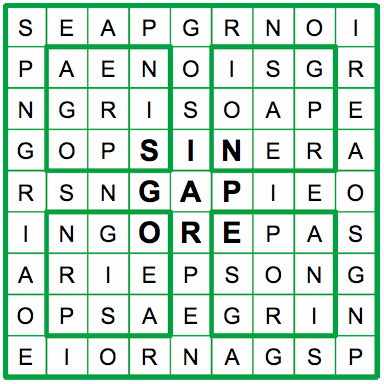}
\par\end{centering}

\protect\caption{\label{fig:Win puzzle 1-9}Win One}
\end{figure}

\noindent \textbf{Win One} (2014\footnote{Year of first publication (see \cite{Latin Puzzles site}).}).
``Complete the board in Fig. \ref{fig:Win puzzle 1-9} so that every
row, column and highlighted region contains numbers $1$ to $9$''.\\

\begin{figure}
\begin{centering}
\includegraphics[height=4.5cm]{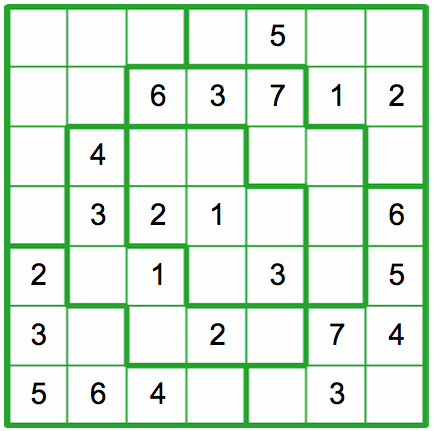}~~~~~~\includegraphics[height=4.5cm]{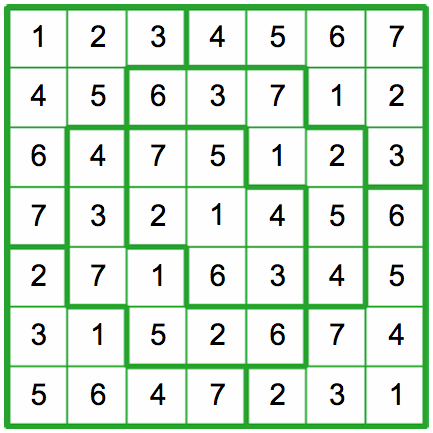}
\par\end{centering}

\protect\caption{\label{fig:Wave One puzzle}Wave One}
\end{figure}

\noindent \textbf{Wave One} (2014). ``Complete the board in Fig.
\ref{fig:Wave One puzzle} so that every row, column and highlighted
region contains numbers $1$ to $7$''. In contrast to \emph{Win
One} puzzles, the extra asterisms here partition the set of cells.
The solution in this case is a type of Latin square called \emph{Gerechte
square} \cite{Gerechte squares}: a regular Latin square with an extra
partition of the board in $n$ regions (seven here) with $n$ cells
each.\\

\begin{figure}
\begin{centering}
\includegraphics[height=4.5cm]{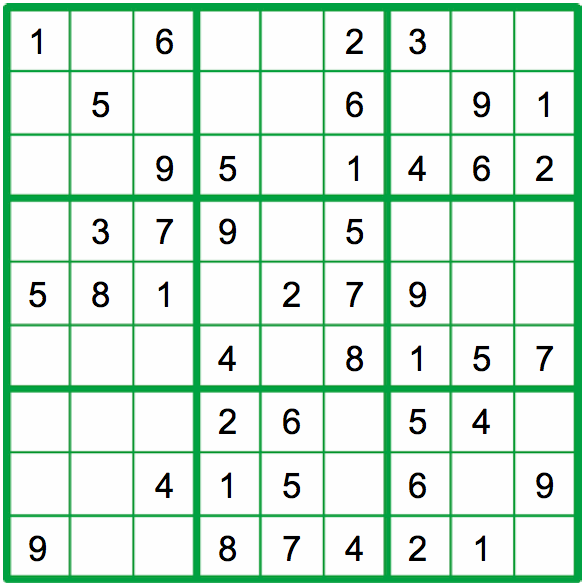}~~~~~~\includegraphics[height=4.5cm]{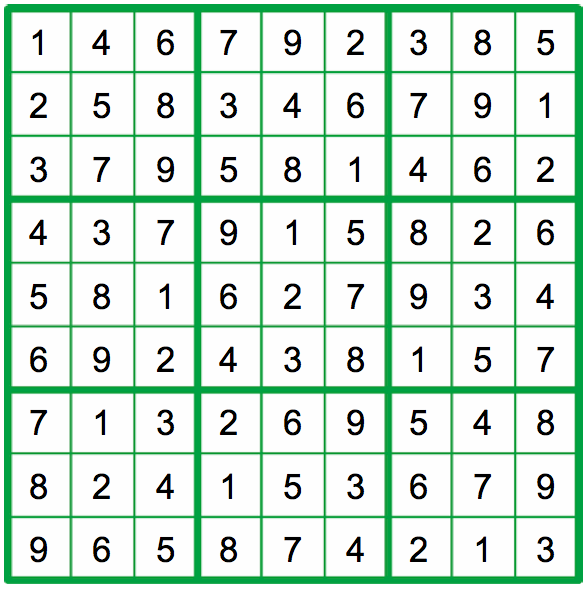}
\par\end{centering}

\protect\caption{\label{fig:Sudoku puzzle-1}Sudoku}
\end{figure}

\noindent \textbf{Sudoku} (1979, see \cite{Howard Garns}). ``Complete
the board in Fig. \ref{fig:Sudoku puzzle-1} so that every row, column
and highlighted region contains numbers $1$ to $9$''\footnote{\emph{Sudoku} puzzles are then Gerechte squares too. Asterisms in
\emph{Sudoku} (rows, columns and sub-squares) are called \emph{units}
or \emph{scopes.}}. \emph{Sudoku} is arguably the most popular Latin square Puzzle variant,
and is itself the base of other related Puzzles \cite{Taking Sudoku Seriously}.\\

\begin{figure}
\begin{centering}
\includegraphics[height=4.5cm]{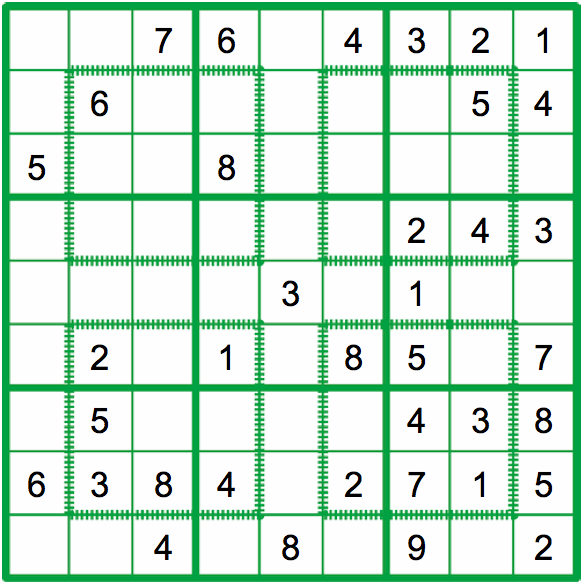}~~~~~~\includegraphics[height=4.5cm]{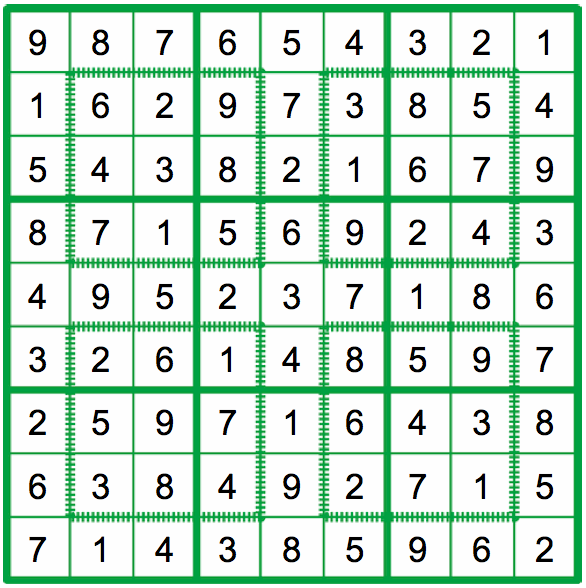}
\par\end{centering}

\protect\caption{\label{fig:Quadoku puzzle 1-9}Quadoku One}
\end{figure}

\noindent \textbf{Quadoku One} (2015). ``Complete the board in Fig.
\ref{fig:Quadoku puzzle 1-9} so that every row, column and highlighted
region contains numbers $1$ to $9$''. This Puzzle combines \textit{Sudoku}
and \textit{Win One} Puzzles into one\footnote{\noindent The \textit{Quadoku One} Puzzle is also called \emph{Hyperdoku}.}.

\subsection{Board Shape}

The next generalization of the Latin square Puzzle is just aesthetic:
the board has a non-square shape, but is topologically equivalent
to a Latin square.\\

\begin{figure}
\begin{centering}
\includegraphics[height=5cm]{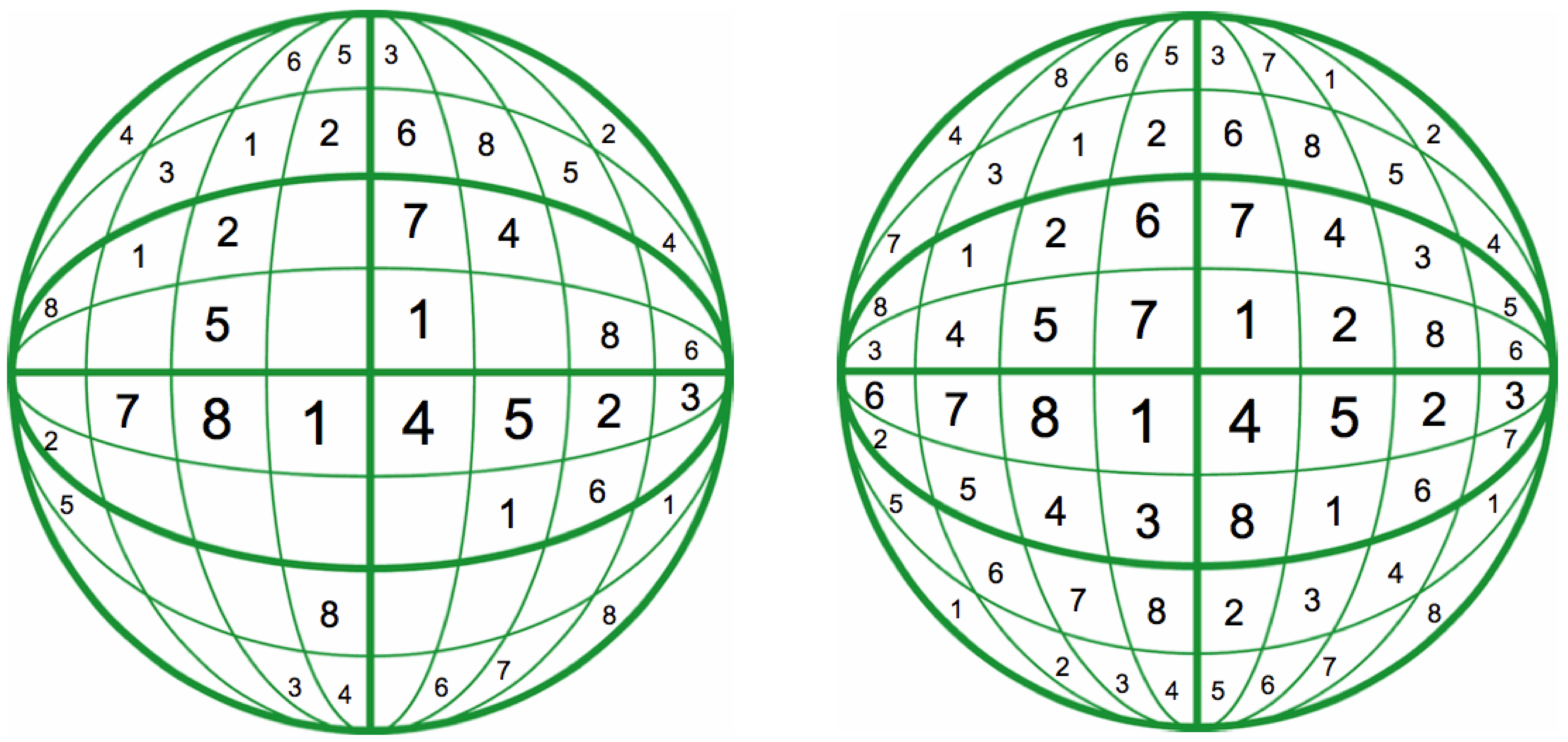}
\par\end{centering}

\protect\caption{\label{fig:Orb One puzzle}Orb One}
\end{figure}

\noindent \textbf{Orb One} (2014). ``Complete the board in Fig. \ref{fig:Orb One puzzle}
so that every parallel, meridian and  highlighted region contains
numbers $1$ to $8$''. The solution here is topologically equivalent
also to a Gerechte square.\\

\subsection{Board Symmetry\label{sub:Board-Symmetry}}

This generalization opens the door to Puzzles whose solutions are
not Latin squares.

\subsubsection{Latin Polygons}

\paragraph{\noindent Dihedral Groups}

\noindent Any of the following eight geometric transformations leaves
a square invariant:
\begin{itemize}
\item 0º, 90º, 180º and 270º rotations
\item 2 reflections across the two diagonals
\item 2 reflections, one across the vertical symmetry axis, another across
the horizontal one
\end{itemize}
Together with the operation of composition of transformations they
form the so-called \emph{dihedral group} \textit{$D_{4}$} \cite{Regular polytopes Coxeter},
the structure that captures the symmetry of the square. Similarly,
every regular polygon of $n$ sides has an associated $D_{n}$ dihedral
group with $n$ rotations and $n$ reflections.

When a $D_{4}$ transformation is applied to a Latin square, another
\textendash not necessarily different\textendash{} Latin square results.
Similarly, a Latin board is a \emph{Latin polygon} if every element
of a particular dihedral group transforms it into another Latin board
\cite{Latin Polytopes arxiv}. This is the case with the puzzle solutions
that follow (see \cite{Latin Polytopes arxiv} for their specific
definitions).\\

\begin{figure}
\begin{centering}
\includegraphics[height=5cm]{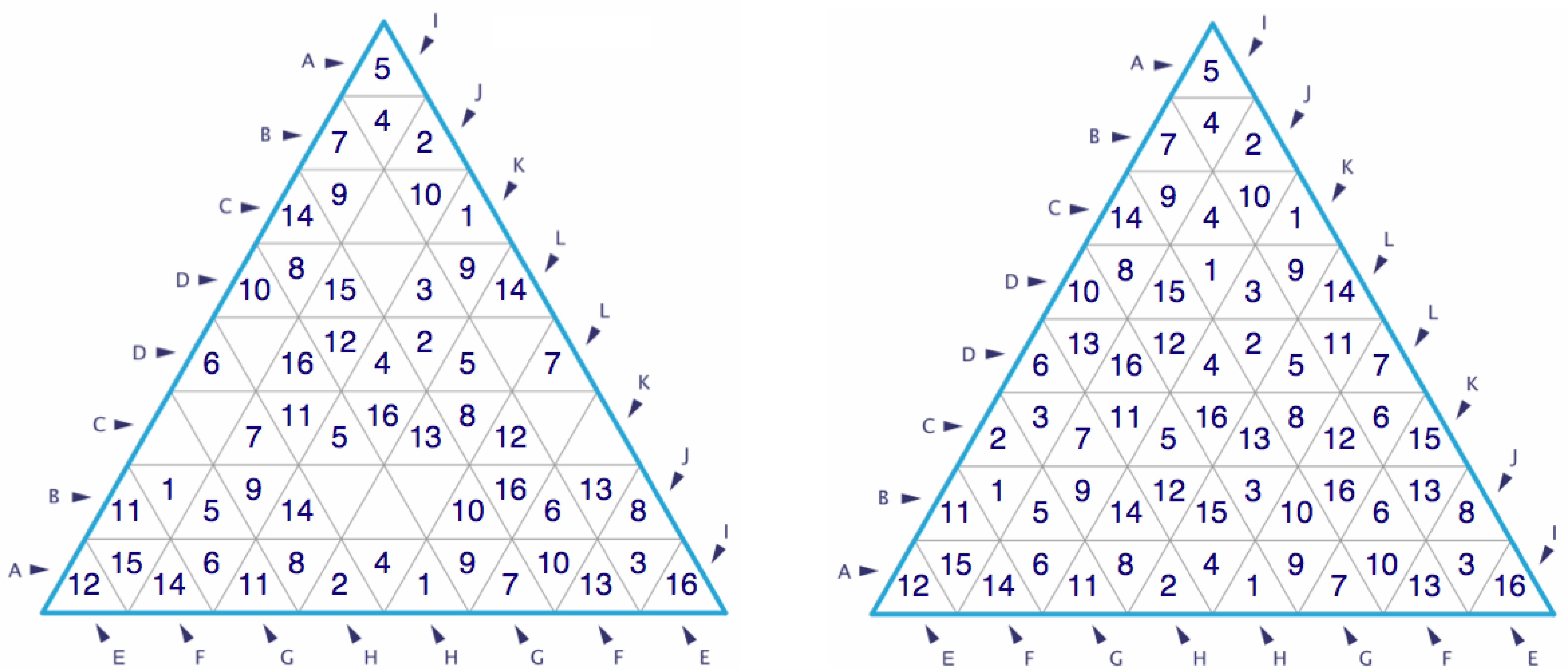}
\par\end{centering}

\protect\caption{\label{fig:Canario One puzzle}Canario One}
\end{figure}

\noindent \textbf{Canario One} (2012). ``Complete the board in Fig.
\ref{fig:Canario One puzzle} so that numbers 1 to 16 appear on every
pair of triangle stripes pointed to by the same letter''. Solutions
to \textit{Canario} puzzles are \emph{Latin triangles}.\\

\begin{figure}
\begin{centering}
\includegraphics[height=5cm]{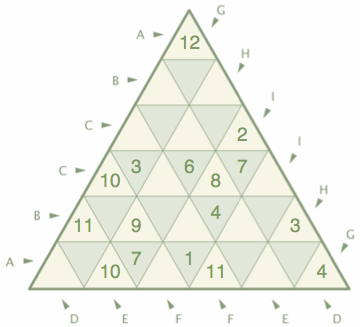}~~~\includegraphics[height=5cm]{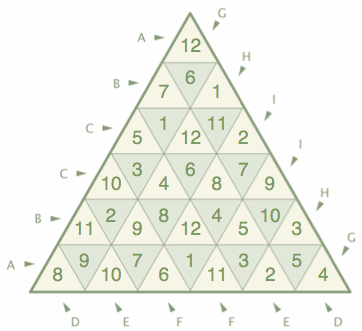}
\par\end{centering}

\protect\caption{\label{fig:Monthai One puzzle}Monthai One}
\end{figure}

\noindent \textbf{Monthai One} (2013). ``Complete the board in Fig.
\ref{fig:Monthai One puzzle} so that numbers 1 to 12 appear on every
pair of triangle stripes pointed to by the same letter''. Solutions
to \textit{Monthai} puzzles are also Latin triangles.\\

\begin{figure}
\begin{centering}
\includegraphics[height=5cm]{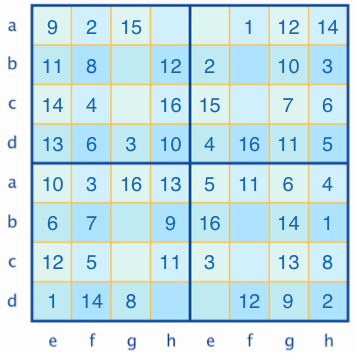}~~~~~~\includegraphics[height=5cm]{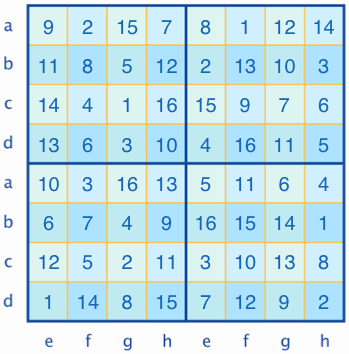}
\par\end{centering}

\protect\caption{\label{fig:Tartan One puzzle}Tartan One}
\end{figure}

\noindent \textbf{Tartan One} (2013). ``Complete the board in Fig.
\ref{fig:Tartan One puzzle} so that numbers 1 to 16 appear on every
pair of rows pointed to by the same letter, every pair of columns
pointed to by the same letter and every highlighted sub-square''.
Solutions to \emph{Tartan} puzzles are \emph{free Latin squares}.\\

\begin{figure}
\begin{centering}
\includegraphics[height=5.8cm]{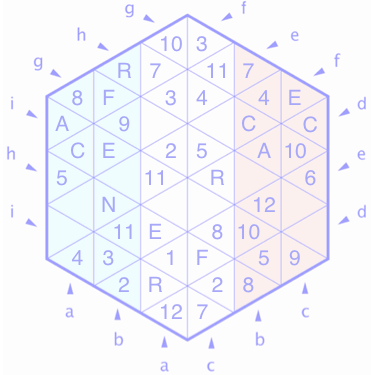}~~~~\includegraphics[height=5.8cm]{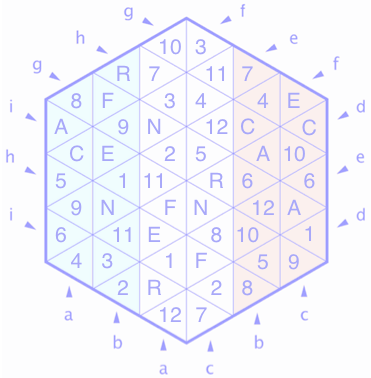}
\par\end{centering}

\protect\caption{\label{fig:Douze France One puzzle}Douze France One}
\end{figure}

\noindent \textbf{Douze France One} (2013). ``Complete the board
in Fig. \ref{fig:Douze France One puzzle} so that the numbers 1 to
12 and the letters F, R, A, N, C, E appear on every pair of triangle
stripes pointed to by the same letter''. Solutions to \textit{Douze
France} puzzles are \emph{Latin hexagons}.\\

\begin{figure}
\begin{centering}
\includegraphics[height=6cm]{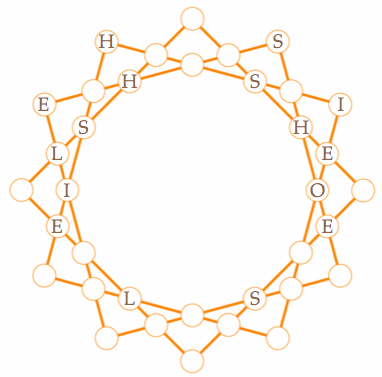}~~~~~\includegraphics[height=6cm]{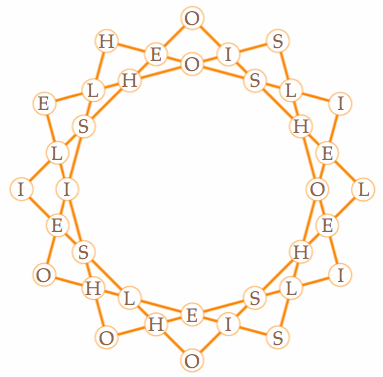}
\par\end{centering}

\protect\caption{\label{fig:Helios One puzzle}Helios One}
\end{figure}

\noindent \textbf{Helios One} (2013). ``Complete the board in Fig.
\ref{fig:Helios One puzzle} so that the letters H, E, L, I, O, S
appear on every line. Solutions to \textit{Helios} puzzles are \emph{Latin
dodecagons}.

\subsubsection{Latin Polyhedra\label{sub:Latin-polyhedra}}

\paragraph{Polyhedral Groups}

Every regular polyhedron has an associated symmetry group whose elements
leave the polyhedron invariant \cite{Regular polytopes Coxeter}.
Table \ref{tab:symmetryGroupsPlatonicSolids} lists the so called
\emph{polyhedral groups} and the corresponding five platonic solids.
A Latin board is a \emph{Latin polyhedron} if every element of a particular
polyhedral group transforms it into another Latin board \cite{Latin Polytopes arxiv}.
This is the case with the puzzle solutions that follow.

Once folded along the inner black lines and glued along the outer
ones, each of the pictures that follows becomes a polyhedral board
with points and asterisms on its surface, each asterism being delimited
by two lines of the same color. The corresponding puzzles are critical
(see Def. \ref{def:critical-Latin-puzzle}), rated \emph{difficult}
(see Sect. \ref{sub:Rating fair Latin puzzles}) and have a Latin
polyhedron as a solution. \medskip{}
\begin{table}
\begin{centering}
\begin{tabular}{ccc}
\noalign{\vskip0.3cm}
\textbf{\footnotesize{}Group} & \textbf{\footnotesize{}Regular Polyhedron} & \textbf{\footnotesize{}Order}\tabularnewline
\hline 
\noalign{\vskip0.1cm}
{\footnotesize{}tetrahedral} & {\footnotesize{}tetrahedron} & {\footnotesize{}24}\tabularnewline
\noalign{\vskip0.1cm}
{\footnotesize{}octahedral} & {\footnotesize{}cube, octahedron} & {\footnotesize{}48}\tabularnewline
\noalign{\vskip0.1cm}
{\footnotesize{}icosahedral} & {\footnotesize{}dodecahedron, icosahedron} & {\footnotesize{}120}\tabularnewline
\end{tabular}
\par\end{centering}

\protect\caption{\label{tab:symmetryGroupsPlatonicSolids}Polyhedral symmetry groups
of the five Platonic solids and their order}
\end{table}

\begin{figure}
\begin{centering}
\includegraphics[height=5cm]{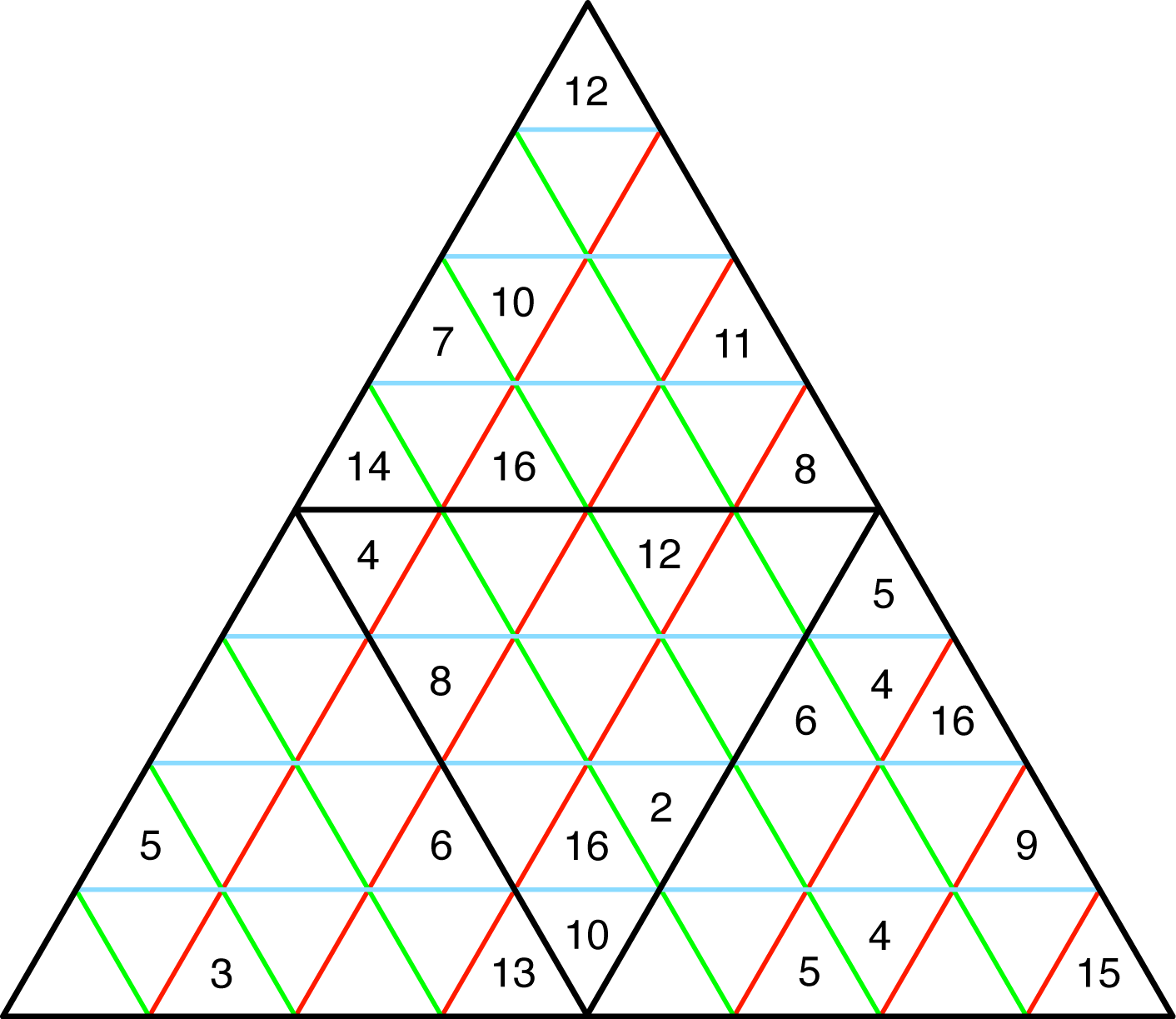}~~~~~~\includegraphics[height=5cm]{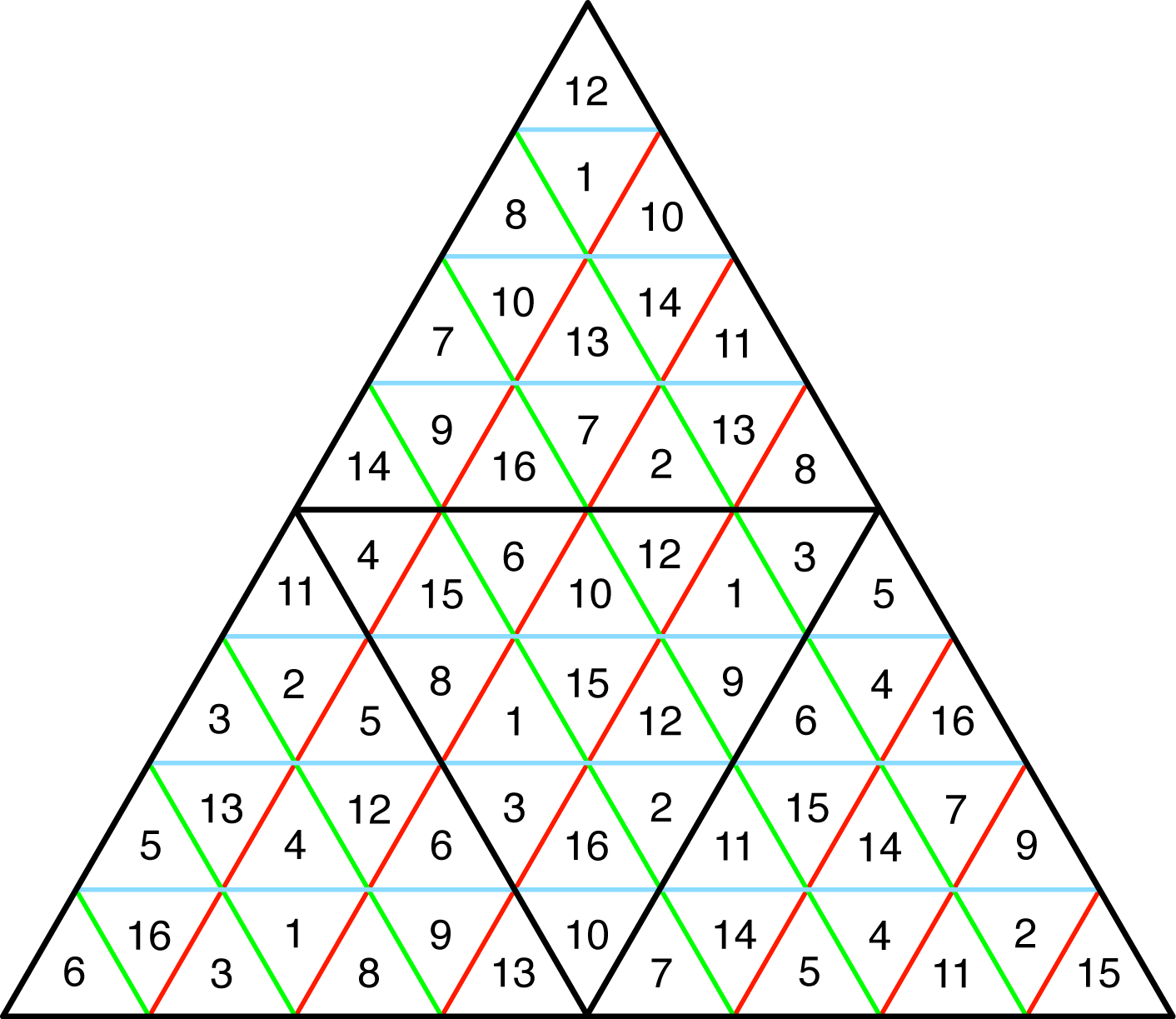}
\par\end{centering}

\protect\caption{\label{fig:Tetra One puzzle}Latin Tetrahedron}

\end{figure}

\noindent \textbf{Latin Tetrahedron} \textbf{Puzzle} (2012). Fig.
\ref{fig:Tetra One puzzle} shows a tetrahedron with twenty-four closed
stripes of triangles on its surface. The challenge is to complete
it so that every stripe holds all numbers from $1$ to $16$.

\noindent 

\begin{figure}
\begin{centering}
\includegraphics[height=8cm]{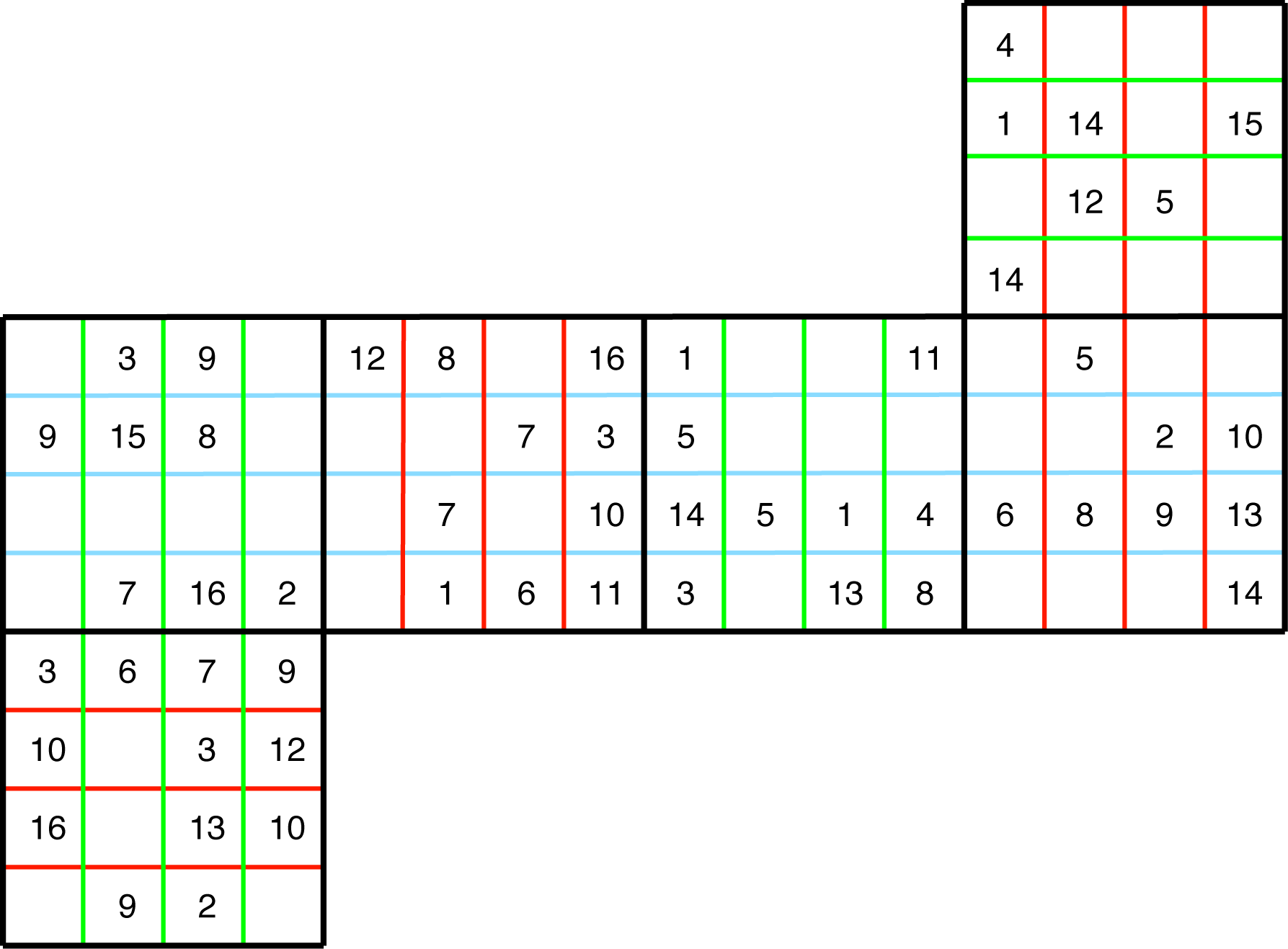}
\par\end{centering}

\begin{centering}
\bigskip{}

\par\end{centering}

\begin{centering}
\bigskip{}

\par\end{centering}

\begin{centering}
\includegraphics[height=8cm]{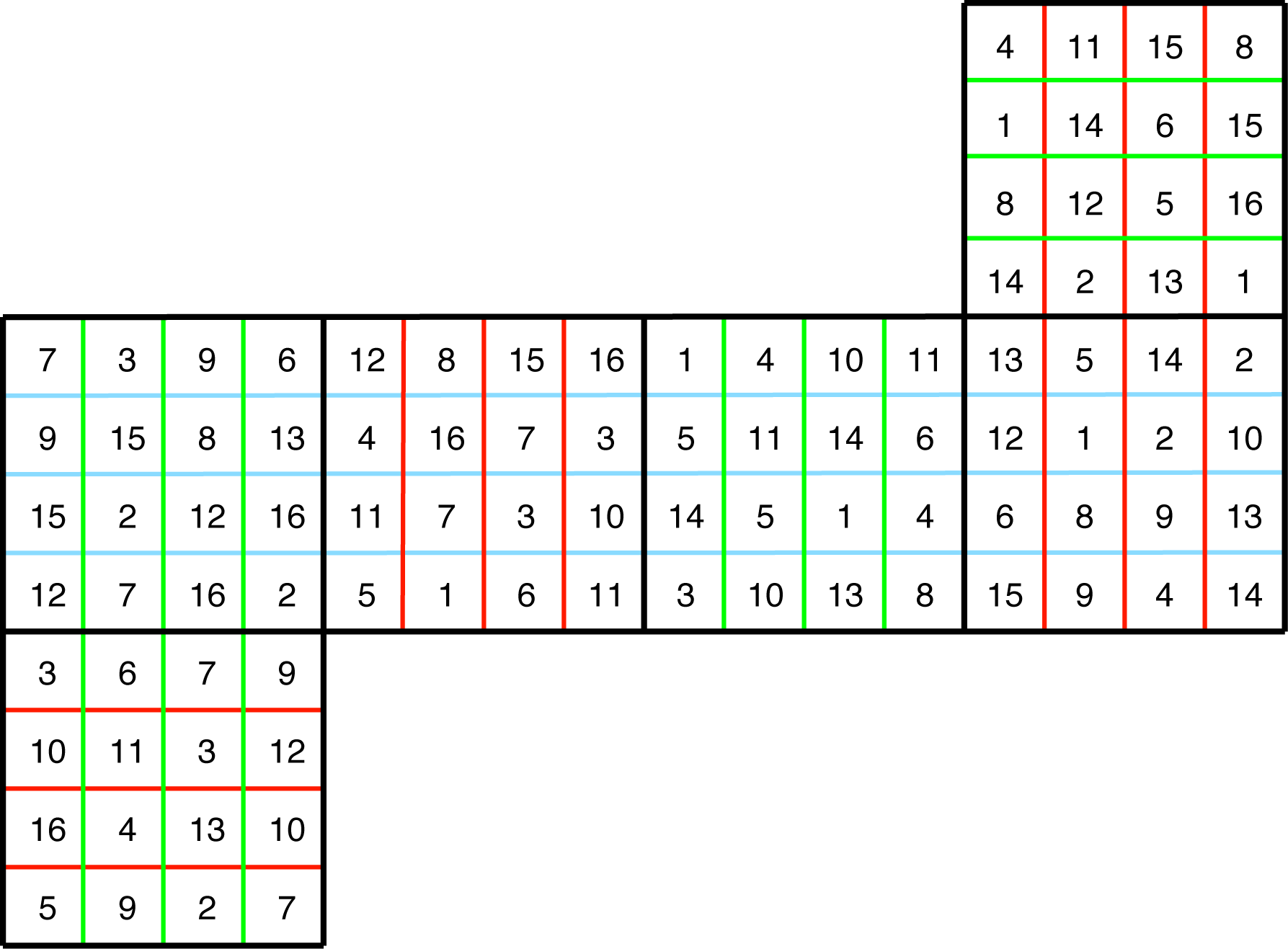}
\par\end{centering}

\protect\caption{\label{fig:Cube One puzzle}Latin Cube}

\end{figure}

\noindent \textbf{Free Latin Cube} \textbf{Puzzle} (2014). Fig. \ref{fig:Cube One puzzle}
shows a cube with twelve closed stripes of squares on its surface.
The challenge is to complete it so that every stripe holds all numbers
from $1$ to $16$.\\

\begin{figure}
\begin{centering}
\includegraphics[height=8cm]{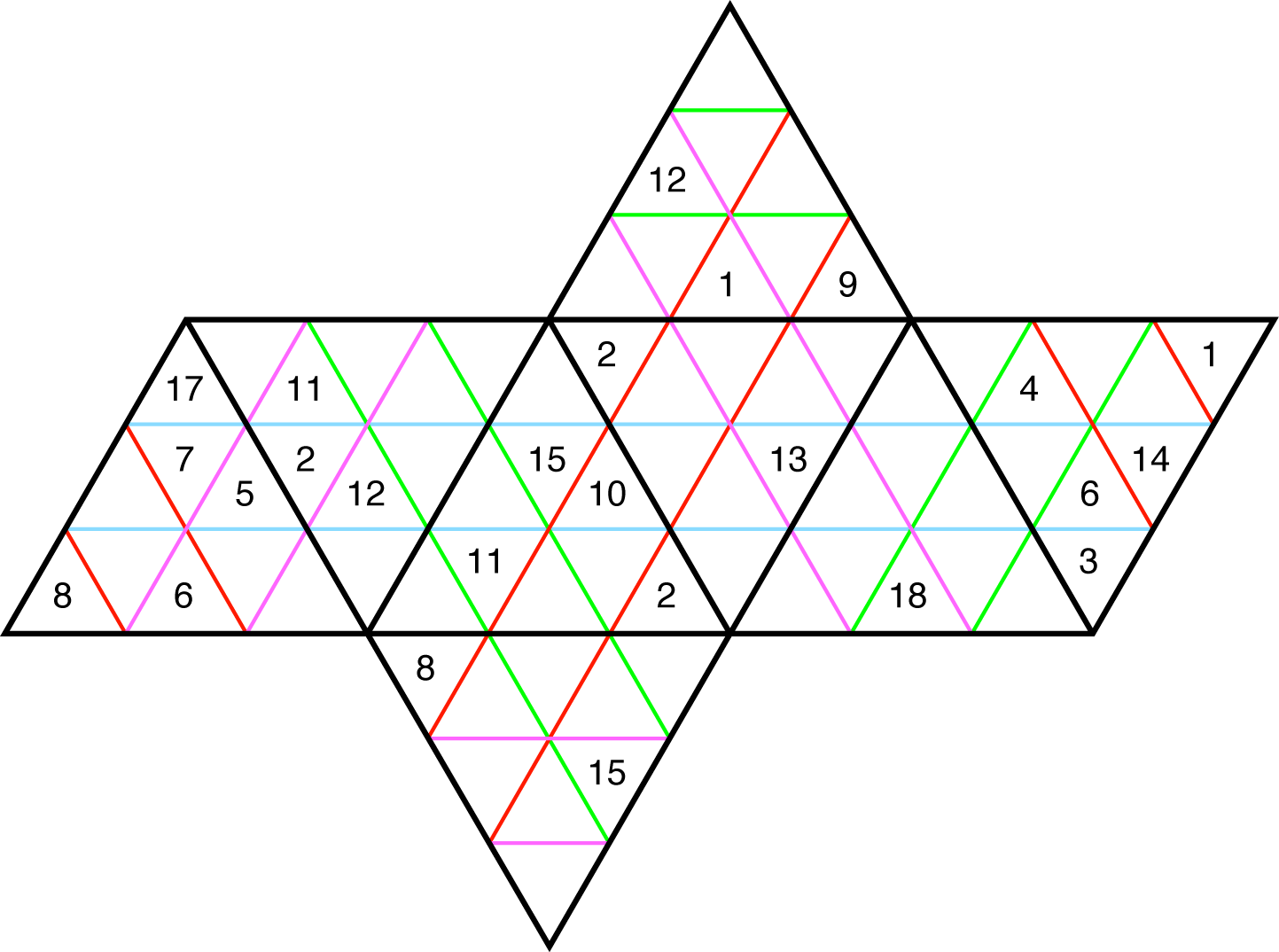}
\par\end{centering}

\begin{centering}
\bigskip{}

\par\end{centering}

\begin{centering}
\bigskip{}

\par\end{centering}

\begin{centering}
\includegraphics[height=8cm]{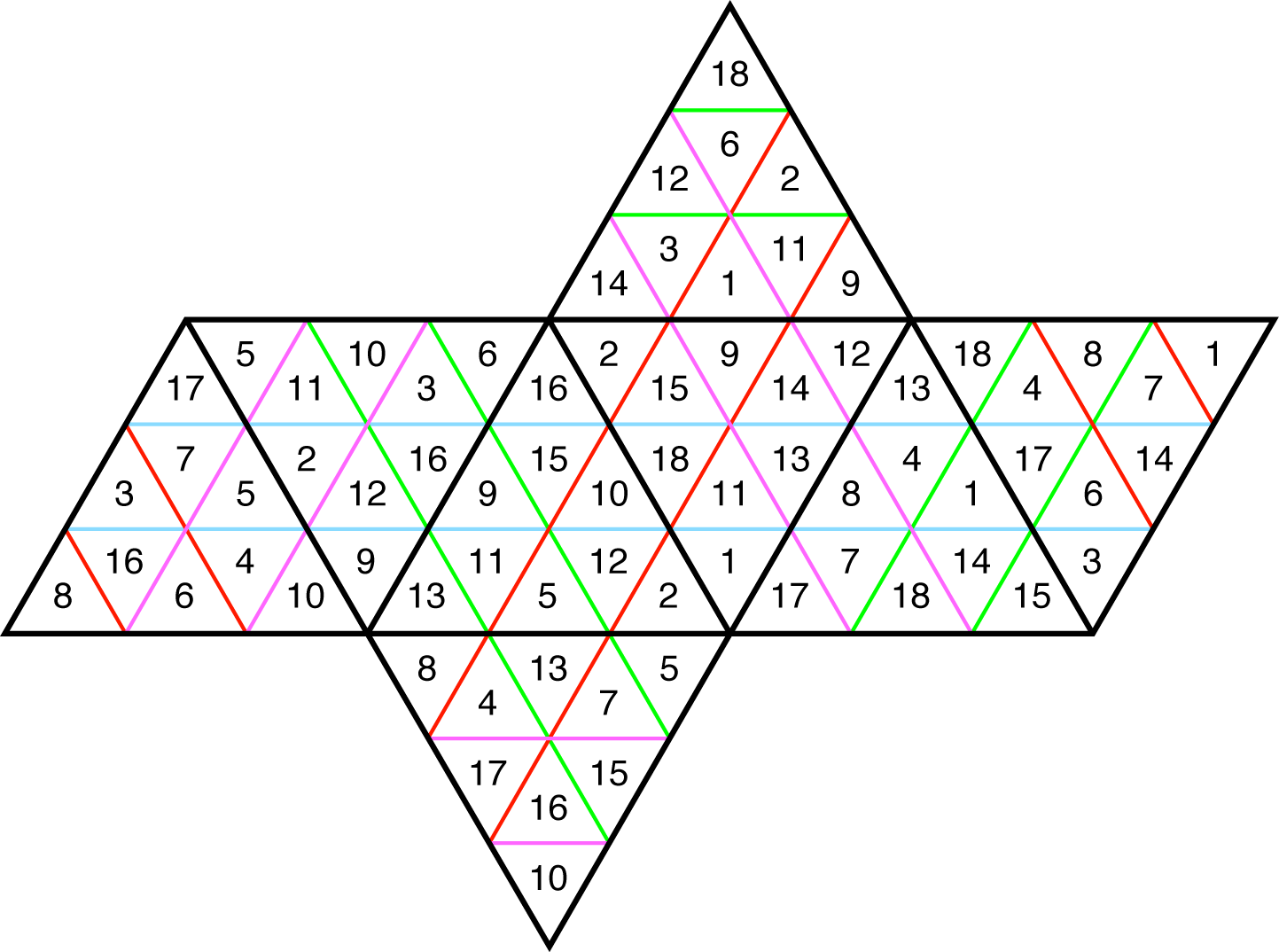}
\par\end{centering}

\protect\caption{\label{fig:Octa One puzzle}Latin Octahedron}

\end{figure}

\noindent \textbf{Latin Octahedron} \textbf{Puzzle} (2014). Fig. \ref{fig:Octa One puzzle}
shows an octahedron with twelve closed stripes of triangles on its
surface. The challenge is to complete it so that every stripe holds
all numbers from $1$ to $18$.\\

\begin{figure}
\begin{centering}
\includegraphics[height=6cm]{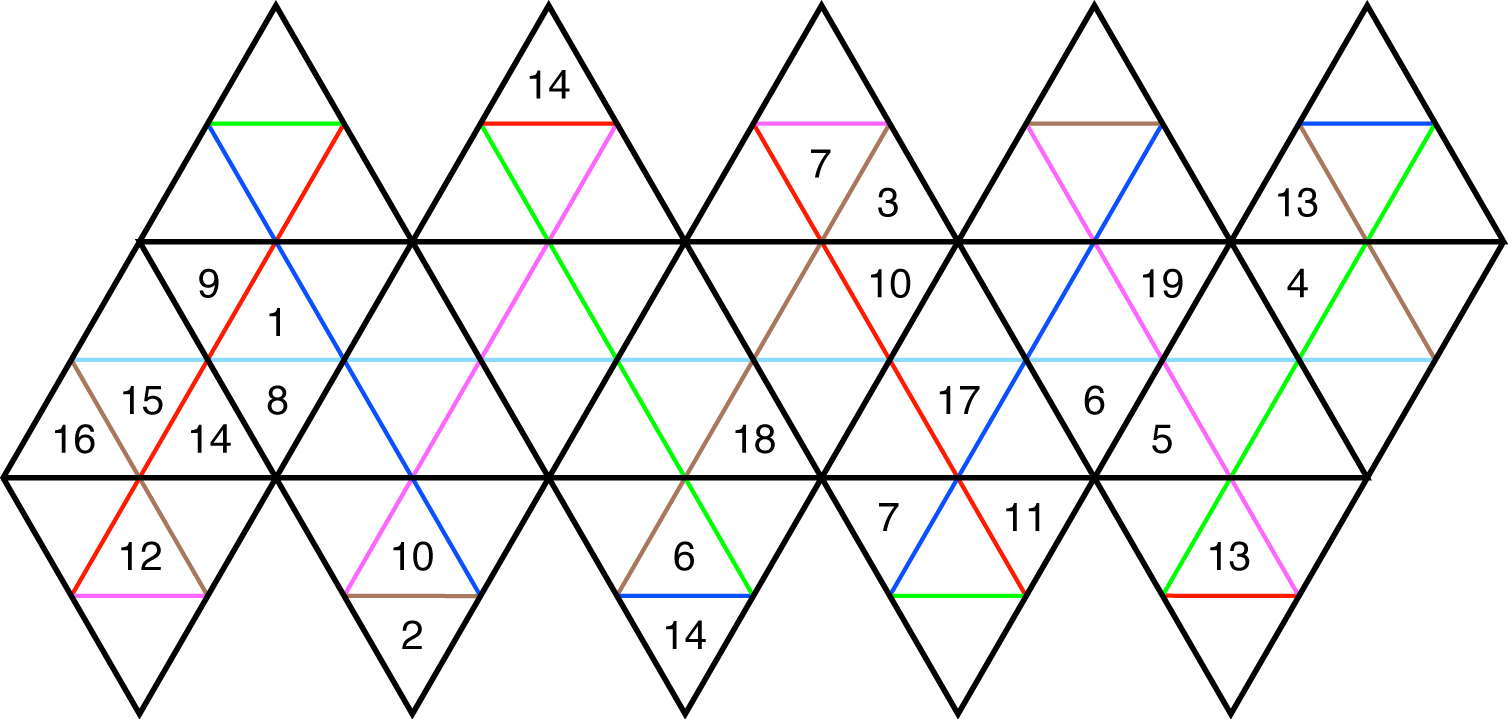}
\par\end{centering}

\begin{centering}
\bigskip{}

\par\end{centering}

\begin{centering}
\bigskip{}

\par\end{centering}

\begin{centering}
\bigskip{}

\par\end{centering}

\begin{centering}
\bigskip{}

\par\end{centering}

\begin{centering}
\includegraphics[height=6cm]{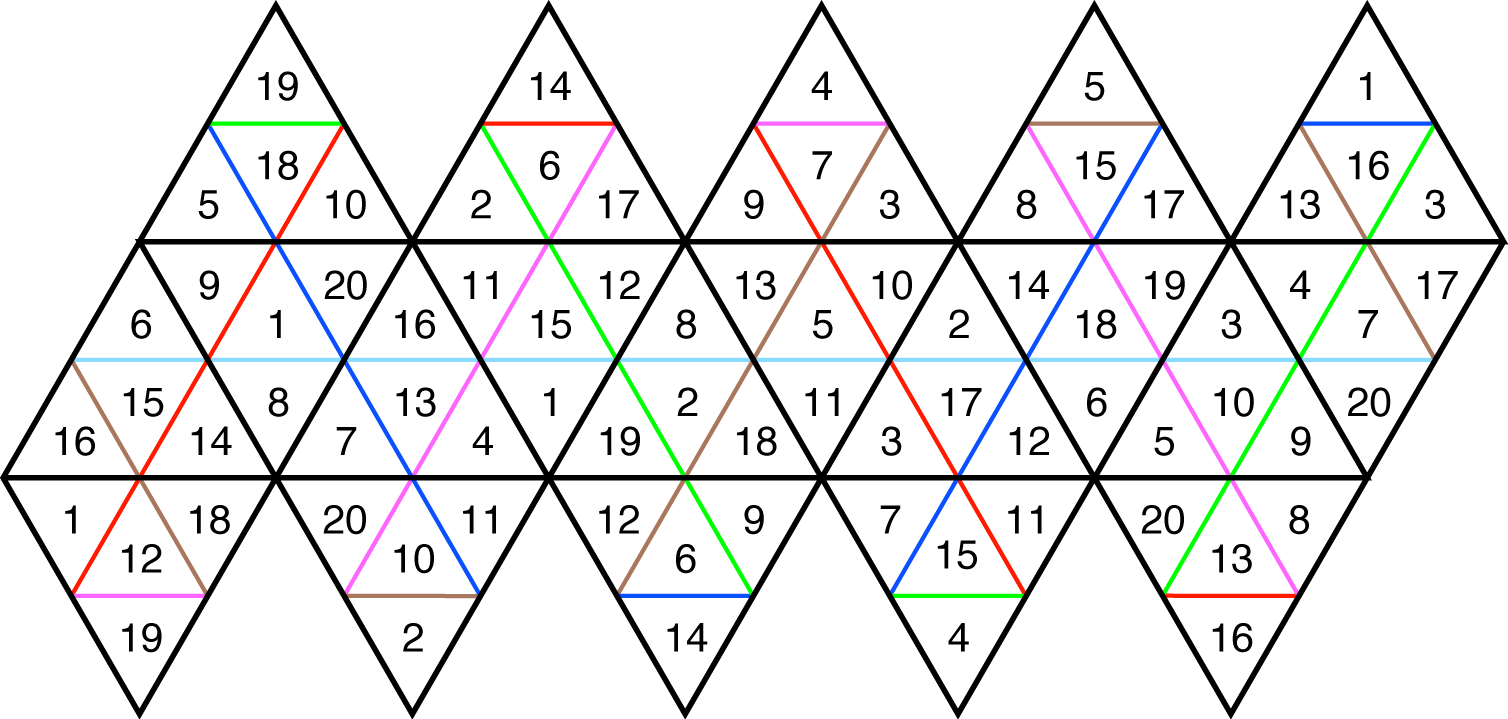}
\par\end{centering}

\protect\caption{\label{fig:Ico One puzzle}Latin Icosahedron}

\end{figure}

\noindent \textbf{Latin Icosahedron} \textbf{Puzzle} (2014). Fig.
\ref{fig:Ico One puzzle} shows an icosahedron with twelve closed
stripes of triangles on its surface. The challenge is to complete
it so that every stripe holds all numbers from $1$ to $20$.\\

\begin{figure}
\begin{centering}
\includegraphics[height=6cm]{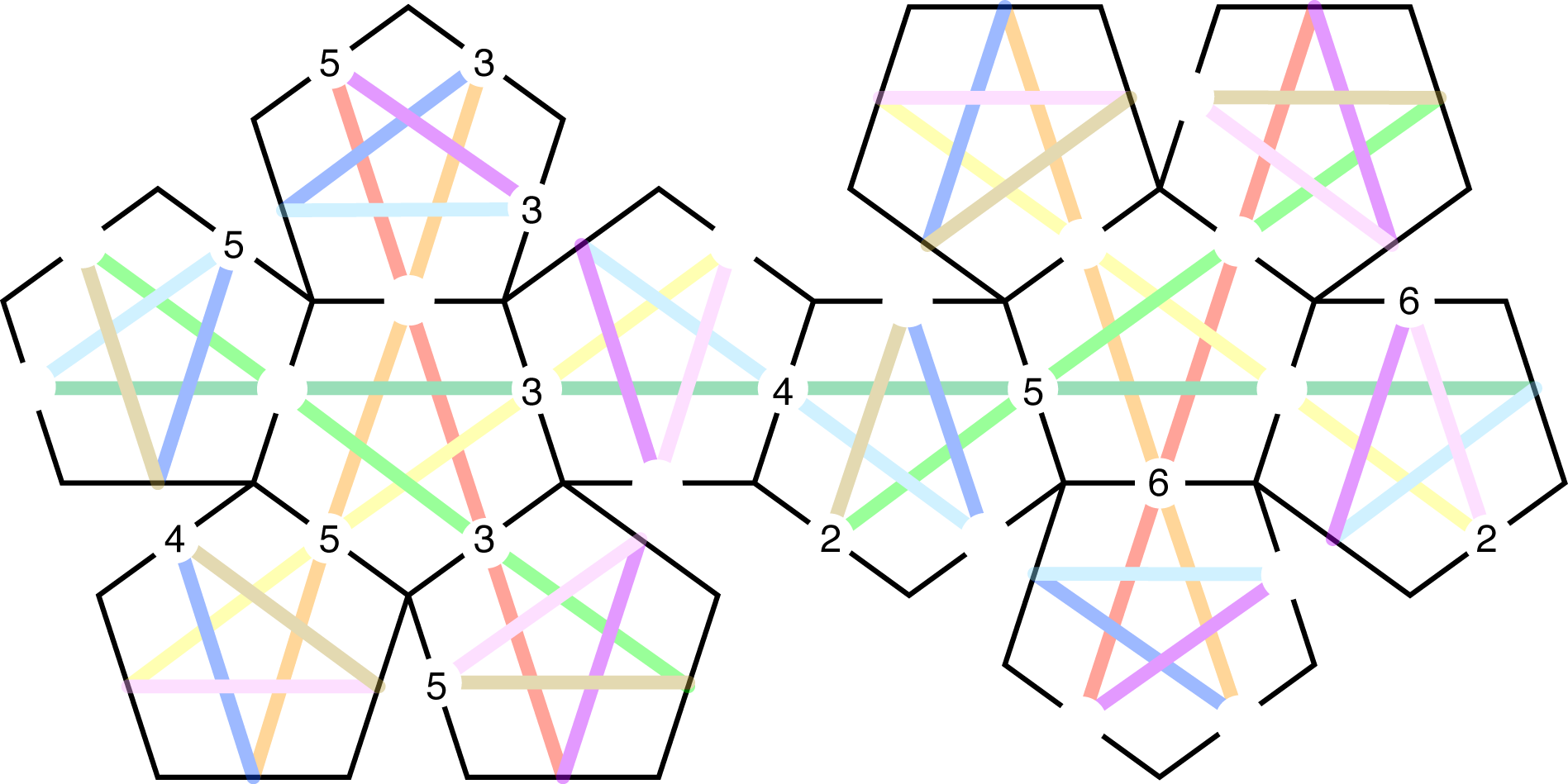}\bigskip{}

\par\end{centering}

\begin{centering}
\bigskip{}

\par\end{centering}

\begin{centering}
\bigskip{}

\par\end{centering}

\begin{centering}
\bigskip{}

\par\end{centering}

\begin{centering}
\includegraphics[height=6cm]{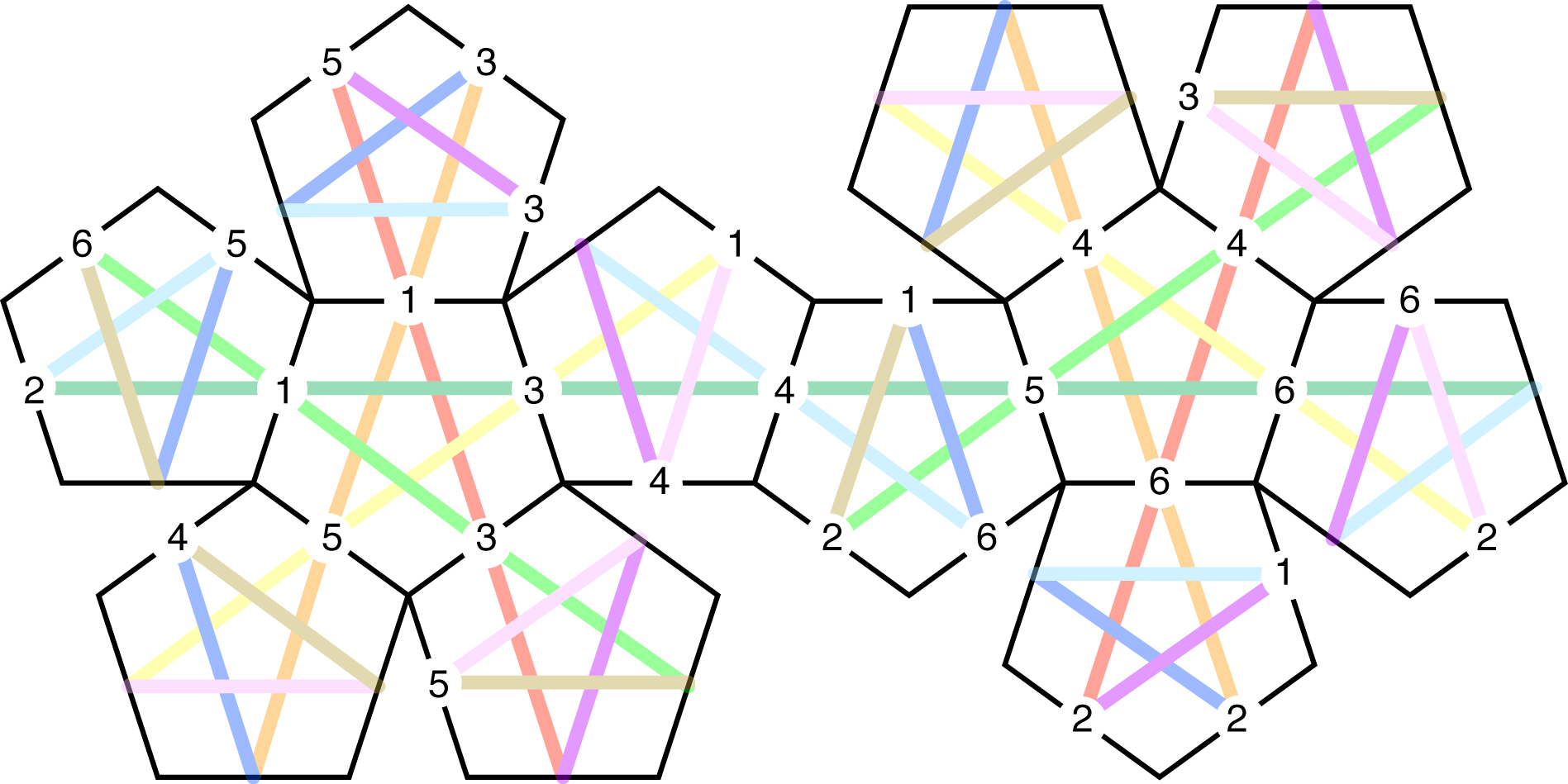}
\par\end{centering}

\protect\caption{\label{fig:Dode One puzzle}Latin Dodecahedron}
\end{figure}

\noindent \textbf{Latin Dodecahedron} \textbf{Puzzle} (2014). Fig.
\ref{fig:Dode One puzzle} shows a dodecahedral board with points
in the middle of every edge and asterisms indicated by lines of the
same color. Every line here must hold all numbers from $1$ to $6$.

\subsection{Repeated Labels\label{sub:Repeated-Labels}}

This generalization allows asterisms to hold repeated labels (see
Def. \ref{def:Latin board definition}).\medskip{}

\begin{figure}
\begin{centering}
\includegraphics[height=4.5cm]{Sudoku_Repeto_Clues}~~~~~~\includegraphics[height=4.5cm]{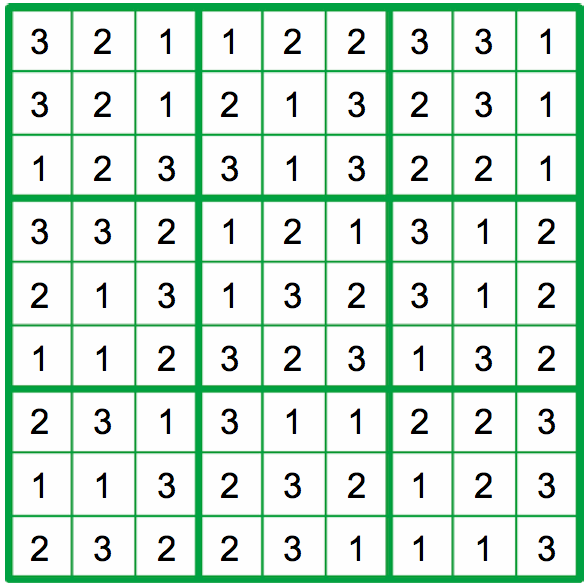}
\par\end{centering}

\protect\caption{\label{fig:Sudoku Ripeto puzzle}Sudoku Ripeto}
\end{figure}

\noindent \textbf{Sudoku Ripeto} (2014). ``Complete the board in
Fig. \ref{fig:Sudoku Ripeto puzzle} so that every row, column and
highlighted region contains the labels $\left\{ 1,1,1,2,2,2,3,3,3\right\} $''.\\

\begin{figure}
\begin{centering}
\includegraphics[height=5cm]{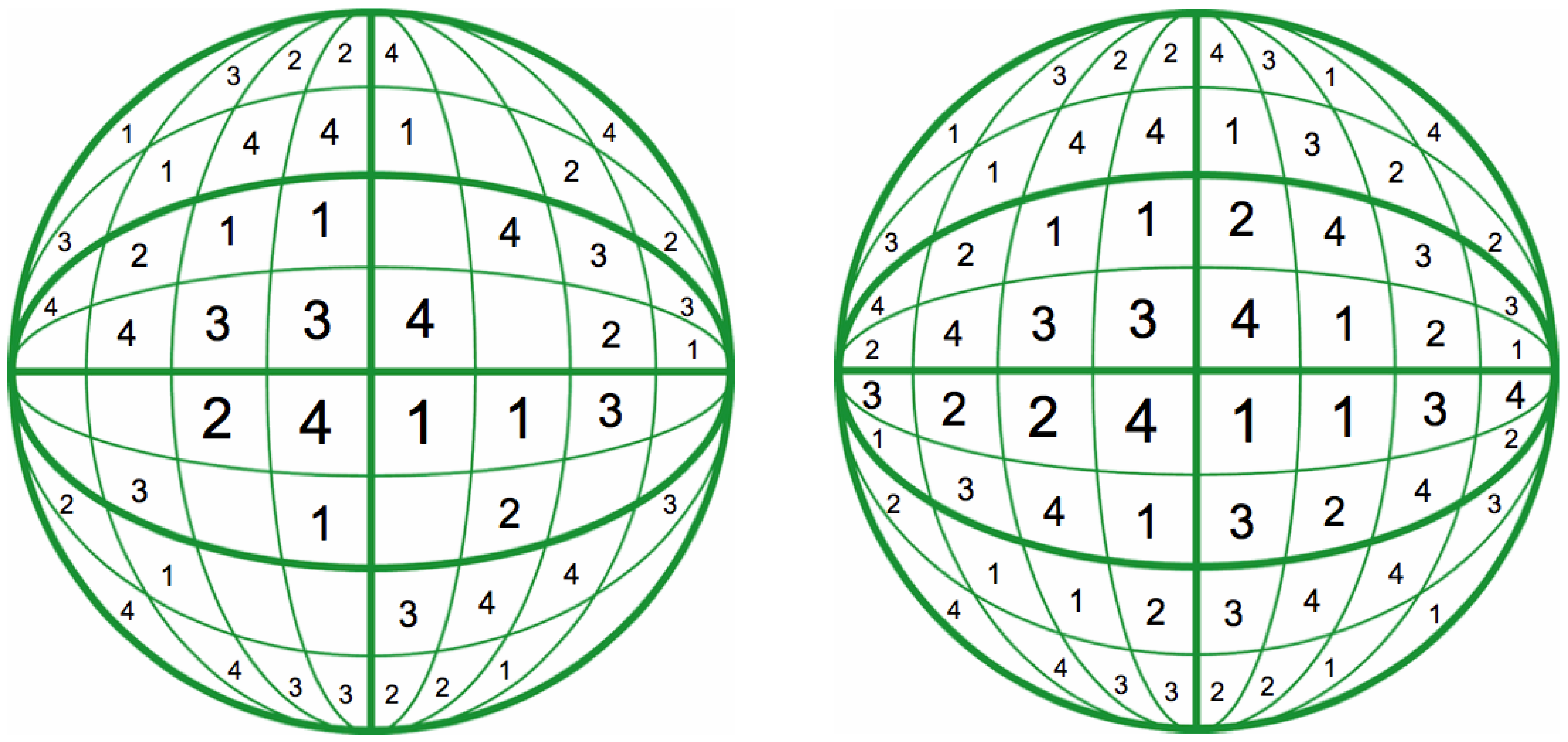}
\par\end{centering}

\protect\caption{\label{fig:Orb Repeto puzzle}Orb Ripeto}
\end{figure}

\noindent \textbf{Orb Ripeto} (2014). ``Complete the board in Fig.
\ref{fig:Orb Repeto puzzle} so that every parallel, meridian and
highlighted region contains the labels $\left\{ 1,1,1,2,2,2,3,3,3\right\} $''.\\

\begin{figure}
\centering{}\includegraphics[height=4.5cm]{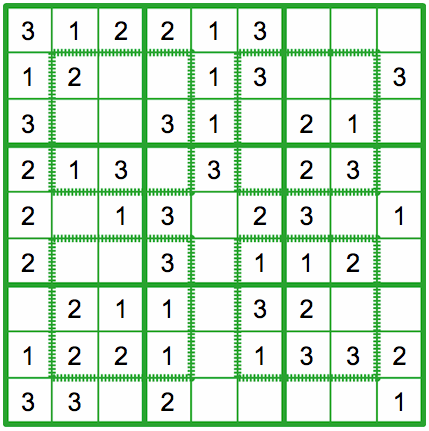}~~~~~~\includegraphics[height=4.5cm]{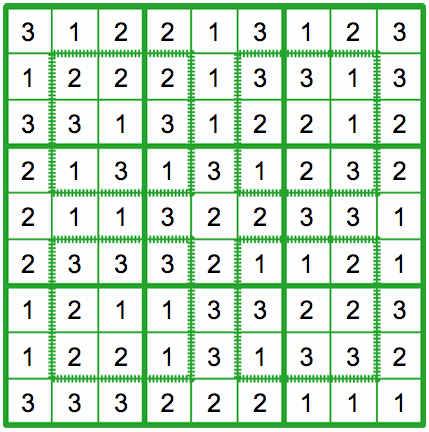}\protect\caption{\label{fig:Quadoku Ripeto puzzle}Quadoku Ripeto}
\end{figure}

\noindent \textbf{Quadoku Ripeto} (2015). ``Complete the board in
Fig. \ref{fig:Quadoku Ripeto puzzle} so that every row, column and
highlighted region contains the labels $\left\{ 1,1,1,2,2,2,3,3,3\right\} $''.

\subsection{Inscriptions\label{sub:Inscriptions}}

Unless otherwise stated, in the examples that follow the challenge
is ``Complete the board in the figure so that every row, column and
highlighted region contains the same letters as the completed row''.
The inscription on the first puzzles is the name of a relevant figure
in the history of \emph{Sudoku} (see Sect. \ref{sec:Introduction}).
\\

\begin{figure}
\begin{centering}
\includegraphics[height=4.5cm]{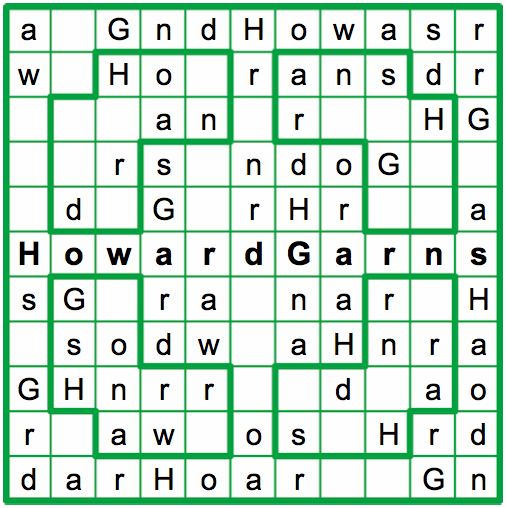}~~~~~~\includegraphics[height=4.5cm]{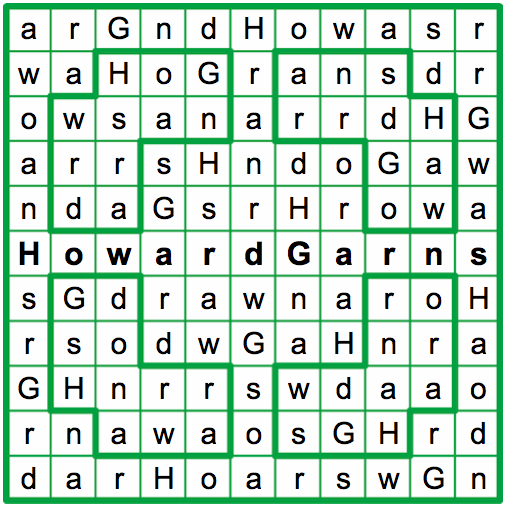}
\par\end{centering}

\protect\caption{\label{fig:Focus puzzle}Custom Focus}
\end{figure}

\noindent \textbf{Custom Focus} (2014). See Fig. \ref{fig:Focus puzzle}.\\

\begin{figure}
\begin{centering}
\includegraphics[height=4.5cm]{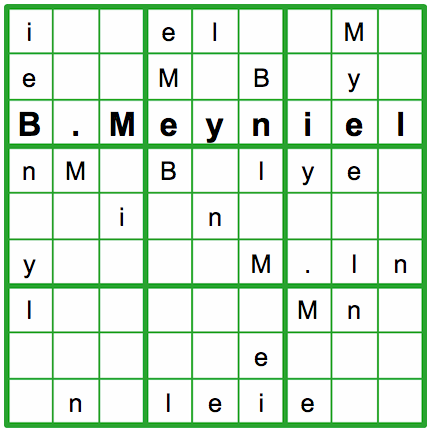}~~~~~~\includegraphics[height=4.5cm]{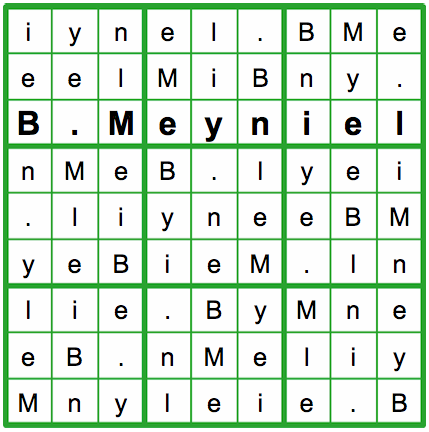}
\par\end{centering}

\protect\caption{\label{fig:Custom Sudoku B.Meyniel}Custom Sudoku}
\end{figure}

\noindent \textbf{Custom Sudoku} (2014). See Fig. \ref{fig:Custom Sudoku B.Meyniel}.\\

\begin{figure}
\begin{centering}
\includegraphics[height=5cm]{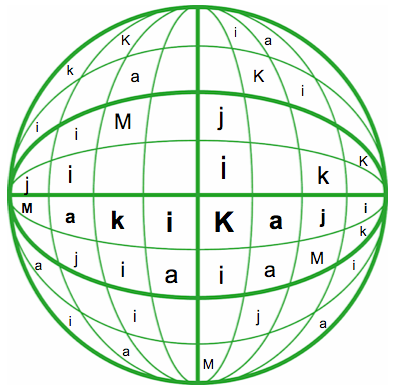}~~~~\includegraphics[height=5cm]{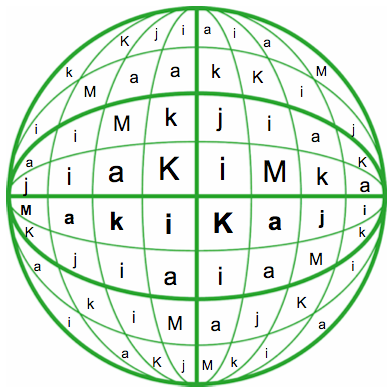}
\par\end{centering}

\protect\caption{\label{fig:Custom Orb Maki Kaji}Custom Orb}
\end{figure}

\noindent \textbf{Custom Orb} (2014). ``Complete the board in Fig.
\ref{fig:Custom Orb Maki Kaji} so that every parallel, meridian and
highlighted region has the same letters as the completed parallel''.\\

\begin{figure}
\begin{centering}
\includegraphics[height=4.5cm]{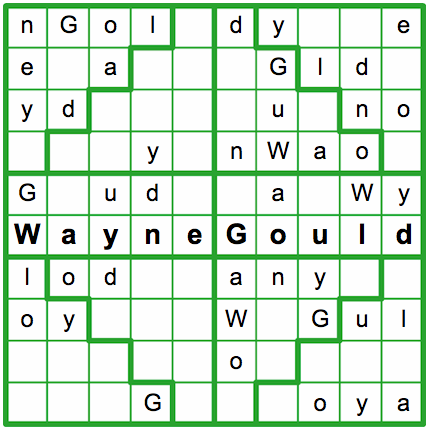}~~~~~~\includegraphics[height=4.5cm]{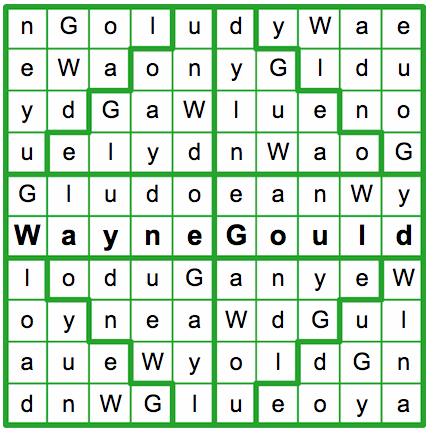}
\par\end{centering}

\protect\caption{\label{fig:CustomSteps Wayne Gould}Custom Steps}
\end{figure}

\noindent \textbf{Custom Steps} (2014). See Fig. \ref{fig:CustomSteps Wayne Gould}.\\

\begin{figure}
\begin{centering}
\includegraphics[height=4.5cm]{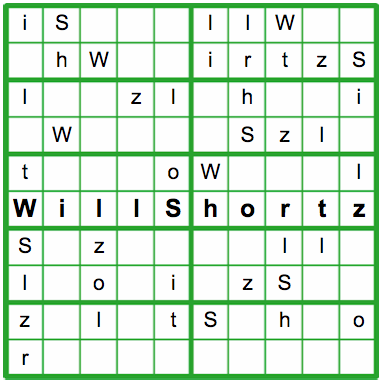}~~~~~~\includegraphics[height=4.5cm]{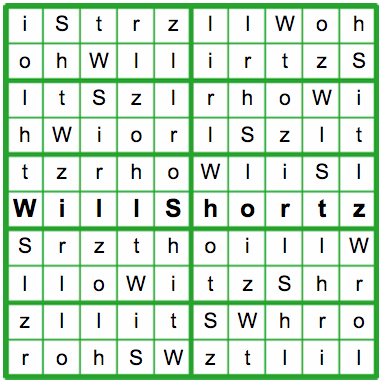}
\par\end{centering}

\protect\caption{\label{fig:Custom Sky Will Shortz}Custom Sky}
\end{figure}

\noindent \textbf{Custom Sky} (2015). See Fig. \ref{fig:Custom Sky Will Shortz}.\\

\begin{figure}
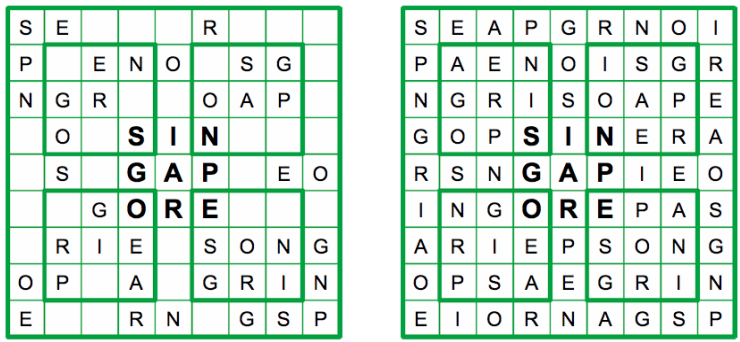

\begin{centering}
\includegraphics[height=4.5cm]{Win_puzzle_clues}~~~~~~\includegraphics[height=4.5cm]{Win_puzzle_solution}
\par\end{centering}

\protect\caption{\label{fig:Custom Win Singapore}Custom Win}
\end{figure}

\noindent \textbf{Custom Win} (2014). ``Complete the board in Fig.
\ref{fig:Custom Win Singapore} so that every row, column and highlighted
region contains the letters in bold typeface''. In this example the
inscription is not confined in an asterism.\\

\begin{figure}
\begin{centering}
\includegraphics[height=4.5cm]{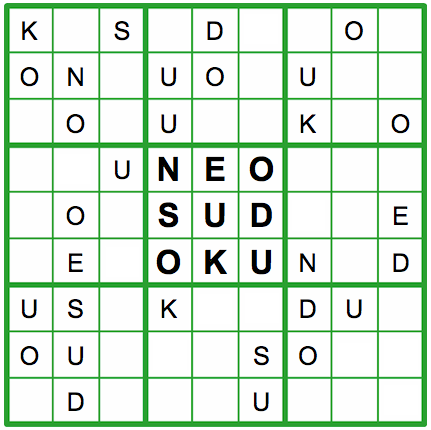}~~~~~~\includegraphics[height=4.5cm]{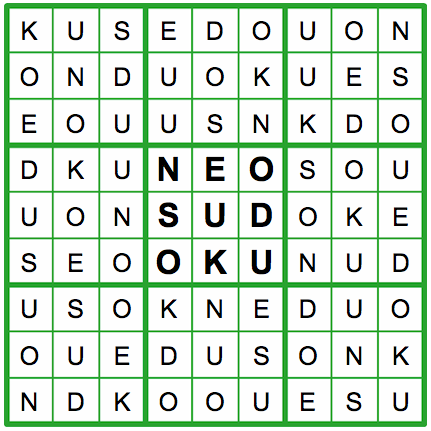}
\par\end{centering}

\protect\caption{\label{fig:didoku NEOSUDOKU-1}Critical Custom Sudoku}
\end{figure}

\noindent \textbf{Custom Sudoku} (2014)\label{Critical Custom-Sudoku}.
``Complete the board in Fig. \ref{fig:didoku NEOSUDOKU-1} so that
every row, column and highlighted region contains the letters in bold
typeface''. This is a critical Latin puzzle (see Def. \ref{def:critical-Latin-puzzle})
with inscription: the removal of any label not in the central sub-square
destroys the puzzle condition. This particular puzzle is rated \emph{difficult}
by the rating method described in Sect. \ref{sub:Rating fair Latin puzzles}.

\section{Finding Latin Boards and puzzles\label{sec:Finding-Latin-puzzles}}

In this section we will use the shorthands indicated in Table \ref{tab:Latin puzzles acronyms}.
References to the term definitions appear in parentheses.

\begin{table}
\begin{raggedright}
\begin{tabular}{rl}
\noalign{\vskip0.1cm}
\textbf{\footnotesize{}Term} & \textbf{\footnotesize{}Acronym}\tabularnewline
\noalign{\vskip0.1cm}
{\footnotesize{}partial labeled board (Sect. \ref{sub:def Partial-Labeled-Board})} & {\footnotesize{}PB}\tabularnewline
\noalign{\vskip0.1cm}
{\footnotesize{}partial Latin board (Sect. \ref{def:partial Latin board})} & {\footnotesize{}PLB}\tabularnewline
\noalign{\vskip0.1cm}
{\footnotesize{}completable partial Latin board (Sect. \ref{sub:Completable-Partial-Latin-Board})} & {\footnotesize{}CPLB}\tabularnewline
\noalign{\vskip0.1cm}
{\footnotesize{}Latin puzzle (Sect. \ref{def:Latin-puzzle})} & {\footnotesize{}LP}\tabularnewline
\noalign{\vskip0.1cm}
{\footnotesize{}Latin board (Sect. \ref{def:Latin board definition})} & {\footnotesize{}LB}\tabularnewline
\end{tabular}
\par\end{raggedright}

\protect\caption{\label{tab:Latin puzzles acronyms}Terms related to Latin boards and
corresponding acronyms}

\end{table}

\subsection{Algorithms}

Algorithm \texttt{list\_LBs\_from\_PB()} in Alg. \ref{pro:LBs from PB}
takes a PB and a positive integer \texttt{N $\ensuremath{\in\mathbb{N}\cup\{\infty\}}$}
as input. Given sufficient memory and time it returns a list with\label{output of list_LBs_from_PB}:
\begin{enumerate}
\item \emph{all} LBs solutions to PB when \texttt{N} $=\infty$
\item \texttt{N }solutions when \texttt{$0<$ N} $<\infty$ and PB has a
total number of solutions $\geq$ \texttt{N}
\item a number of solutions $\leq$ \texttt{N }when \texttt{$0<$ N} $<\infty$
and PB has a total number of solutions $\leq$ \texttt{N}
\end{enumerate}
\noindent If the output list is empty in the first case, then either
the input PB is a PLB but not a CPLB, or it is not a PLB at all. If
it contains a single solution, the input is both a CPLB and an LP.

\begin{algorithm}
\medskip{}
\texttt{\small{}~~algorithm list\_LBs\_from\_PB(PB, N):}{\small \par}

\texttt{\small{}~~~~if PB is not a PLB:}{\small \par}

\texttt{\small{}~~~~~~return empty\_list}{\small \par}

\texttt{\small{}~~~~elif PB is an LB:}{\small \par}

\texttt{\small{}~~~~~~return {[}PB{]}}{\small \par}

\texttt{\small{}~~~~else:}{\small \par}

\texttt{\small{}~~~~~~find LBs for PB as per N and add them
to list\_LBs}{\small \par}

\texttt{\small{}~~~~~~return list\_LBs}{\small \par}

\begin{centering}
\medskip{}

\par\end{centering}

\protect\caption{\label{pro:LBs from PB}Latin boards from a partial labeled board}

\end{algorithm}

Suppose we have another algorithm that, using \texttt{list\_LBs\_from\_PB()}
with \texttt{N} $=\infty$, takes an LB as input and produces a list
of LPs that have that LB as the solution. Alg. \ref{pro:LPs from LB},
that recursively looks for LPs of ever increasing difficulty, is an
example.

\begin{algorithm}
\begin{centering}
\medskip{}

\par\end{centering}

\texttt{\small{}~~algorithm list\_LPs\_from\_LB(LB):						}{\small \par}

\texttt{\small{}~~~~list\_LPs = empty\_list }{\small \par}

\texttt{\small{}~~~~algorithm recur(PLB):							  }{\small \par}

\texttt{\small{}~~~~~~for label in PLB:}{\small \par}

\texttt{\small{}~~~~~~~~remove label from PLB				 }{\small \par}

\texttt{\small{}~~~~~~~~l\_LBs = list\_LBs\_from\_PB(PLB,
N $=\infty$)		 }{\small \par}

\texttt{\small{}~~~~~~~~if l\_LBs contains a single element:   }{\small \par}

\texttt{\small{}~~~~~~~~~~add a copy of PLB to list\_LPs}{\small \par}

\texttt{\small{}~~~~~~~~~~recur(PLB)						}{\small \par}

\texttt{\small{}~~~~~~~~else:}{\small \par}

\texttt{\small{}~~~~~~~~~~restore label in PLB}{\small \par}

\texttt{\small{}~~~~recur(LB)}{\small \par}

\texttt{\small{}~~~~return list\_LPs}{\small \par}

\begin{centering}
\medskip{}

\par\end{centering}

\protect\caption{\label{pro:LPs from LB}Latin puzzles from a Latin board}

\end{algorithm}

Then, given a particular PB, and providing it is a PLB (a $k$-uniform
empty board and a $k$-multiset for example), Alg. \ref{pro:LPs from PB}
will find \texttt{N} or less PBs and their corresponding LPs\footnote{The algorithm could also stop the search when a quantity given as
parameter is reached.}.

\begin{algorithm}
\begin{centering}
\medskip{}

\par\end{centering}

\texttt{\small{}~~algorithm list\_LPs\_from\_PB(PB, N):}{\small \par}

\texttt{\small{}~~~~list\_LBs = list\_LBs\_from\_PB(PB, N)				}{\small \par}

\texttt{\small{}~~~~list\_all\_LPs = empty\_list }{\small \par}

\texttt{\small{}~~~~for LB in list\_LBs:							}{\small \par}

\texttt{\small{}~~~~~~list\_LPs = list\_LPs\_from\_LB(LB)}{\small \par}

\texttt{\small{}~~~~~~list\_all\_LPs += list\_LPs}{\small \par}

\texttt{\small{}~~return list\_all\_LPs}{\small \par}

\begin{centering}
\medskip{}

\par\end{centering}

\protect\caption{\label{pro:LPs from PB}Latin puzzles from a partial labeled board}

\end{algorithm}

Algorithm \texttt{list\_LBs\_from\_PB()} in Alg. \ref{pro:LBs from PB}
assists us in three ways:
\begin{enumerate}
\item to directly find some or all LBs for a PB
\item in particular, with \texttt{N} $=\infty$, to solve an LP given as
input, as this is by definition a PB that is both a PLB and a CPLB
with a single solution
\item with \texttt{N} $=\infty$ and embedded in another algorithm, to find
LPs for an LB
\end{enumerate}
In the first case we usually do not need all solutions. Moreover,
the number of solutions here can be exceedingly large, so we can set
\texttt{N} to a convenient finite value\footnote{Other alternative or complementary option is to limit the algorithm's
running time.}. In the second case, and in contrast with \texttt{N} $=1$, setting
\texttt{N} $=\infty$ additionally allows us to check that the input
is indeed an LP. In the third case the algorithm must run unbridled
and give all solutions.

\subsection{Latin Boards from Partial Labeled Boards\label{sub:Latin-Boards-from-Partial-Boards}}

In this section we discuss ways to implement\texttt{ list\_LBs\_from\_PB()}
in Alg. \ref{pro:LBs from PB}. As seen, finding LBs for a particular
PB is intimately related to finding LPs. While checking whether a
PB is a PLB is easy \textendash and takes linear time in the number
of cells\textendash{} verifying whether a PLB is also a CPLB is hard
in general. For one, verifying whether a generic partial Latin square
is completable is an NP-complete decision problem \cite{Completing partial LS is NP-complete},
so it is unknown whether there is a general, polynomial time algorithm
for the worst cases able to answer this last question. 

Since partial Latin squares are a subset of PLBs, the problem of verifying
PLBs' completability is at least as hard. Consequently, a task equally
difficult is deciding whether or not a PLB is an LP, since LPs are
a subset of CPLBs.
\begin{defn}
The label that an empty cell in a CPLB has in a solution is a \emph{solution
label} for the cell.
\end{defn}

\begin{defn}
Let $x$ be an empty cell in a CPLB with multiset of labels $L$.
Any set $D\left(x\right)\subseteq\left\{ l,l\in L\right\} $ is a
\emph{domain} for $x$ if and only if it includes all its solution
labels.
\end{defn}
\noindent The set $\left\{ l,l\in L\right\} $ is obviously a domain
for any empty cell, called the \emph{initial domain }for the cell.

\begin{defn}
Let $B$ be a CPLB with $n$ empty cells $x_{1},\ldots,x_{n}$. Then
the Cartesian product $\prod_{i=1}^{n}D\left(x_{i}\right)$, where
$D\left(x_{i}\right)$ is a domain for $x_{i}$, is a\emph{ solution
space} of the board.
\end{defn}
\noindent The tuples $\left(l_{1},\ldots,l_{n}\right)$, with $l_{i}\in D\left(x_{i}\right)$,
are called \emph{label vectors, }whereas every $D=\left\{ D\left(x_{i}\right)\right\} $,
the set of all empty cells' domains, is a \emph{board} \emph{domain}.

\begin{defn}
Let $B$ be a CPLB with a set $X$ of empty cells. Let $D$ and $D'$
be two board domains. We say that $D'\subseteq D$ if and only if
$D'\left(x\right)\subseteq D\left(x\right),\forall x\in X$.
\end{defn}
\noindent We say that $D'$ is\emph{ more advanced }or \emph{stronger}
than $D$. Or equivalently, that $D$ is\emph{ less advanced }or \emph{weaker}
than $D'$\label{def:more advanced board domains}. This binary relation
among board domains is clearly a partial order over the set of board
domains.

\subsubsection{Algorithm Specifications\label{sub:Procedure-Specifications}}

In absence of general efficient methods to find LBs, the guiding principles
to design \texttt{list\_LBs\_from\_PB()} could be:
\begin{itemize}
\item accept a parameter for the number of solutions to be obtained
\item verify first that the input is a PLB
\item assume that this PLB is a CPLB
\item perform \emph{attrition}: find more advanced board domains for the
supposed CPLB until solutions are found, either by reducing to singletons
the domains of empty cells or by discarding entire label vectors.
If the domain for an empty cell becomes empty in the process, then
the initial assumption for the PLB being a CPLB was wrong. The algorithm
should return an empty list as soon as this condition is detected
\item optimize: perform attrition in the most efficient way possible, favoring
at each step polynomial-time attrition over exponential attrition,
and less complex polynomial attrition over more complex polynomial
attrition
\item use relevant heuristics to save space, time and complexity
\end{itemize}
When the input to \texttt{list\_LBs\_from\_PB()} is an LP, this approach
has the additional advantage of mimicking the way a human player would
solve the puzzle, an important feature if we want to have compact
written proofs of the solving process\footnote{As we will see later, this is also very useful to both puzzle authors
and players. }.

\subsubsection{Attrition Theorems\label{sub:Attrition-Theorems}}

The reductions performed by \texttt{list\_LBs\_from\_PB()} on board
domains must be safe, in the sense that they must result in other
domains.
\begin{defn}
Let $B$ be a CPLB and $D$ a board domain for it. An \emph{attrition
theorem} $T$ is one that verifies $T(D)\subseteq D$.
\end{defn}
\noindent Attrition theorems never map board domains to less advanced
ones.

\subsubsection{\noindent Attrition Algorithms \label{def:Attrition-Procedures}}

\noindent An \emph{attrition algorithm} is a sequence of operations
that iteratively searches relevant patterns on the board, checks if
each one complies with the conditions to apply an attrition theorem,
then effectively applies it.

\subsubsection{Types of Attrition\label{sub:Types-of-Attrition}}

\begin{figure}
\begin{centering}
\includegraphics[height=4.5cm]{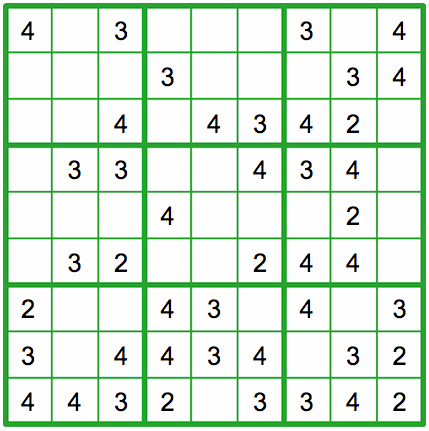}~~~~~~\includegraphics[height=4.5cm]{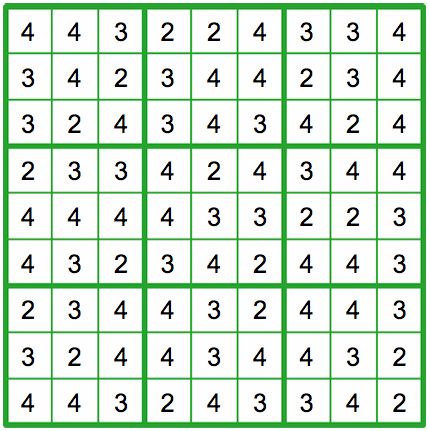}
\par\end{centering}

\protect\caption{\label{fig:sudoku ripeto 223334444}Partial Latin board and a solution}
\end{figure}

The PB in Fig. \ref{fig:sudoku ripeto 223334444} (left) has the same
asterisms as \textit{Sudoku}'s (nine rows, nine columns and nine $3\times3$
sub-squares) and multiset $\left\{ 2,2,3,3,3,4,4,4,4\right\} $. The
board is also a PLB since label counts in every asterism comply with
the multiset. But it is unclear whether it has solutions\footnote{In fact it has only the one shown on its right (this is an instance
of the \emph{Sudoku Ripeto Puzzle} \textendash see Sect. \ref{sub:Repeated-Labels}.}. We study next different attrition methods to find the possible solutions.

\subsubsection{Exponential Attrition by Exhaustive Search}

\noindent With $37$ empty cells, each with initial domain $\left\{ 2,3,4\right\} $,
the initial solution space for the board in Fig. \ref{fig:sudoku ripeto 223334444}
(left) has $3^{37}$ vectors. After the initial checks, one way to
find solutions is to form all possible vectors, write each one on
the board and check compliance, including in the output list all solutions
found along the way. The trivial attrition theorem behind this algorithm
guarantees correctness and termination, but the search performed is
exponential in the number of empty cells. Exhaustive search is the
most clear example of exponential attrition.

\subsubsection{\noindent Exponential Attrition by Backtracking\label{Backtracking}}

\noindent An alternative, although still exponential in nature, is
for \texttt{list\_LBs\_from\_PB()} to implement \emph{backtracking}
\cite{Backtracking}, a technique that systematically checks one,
two, three, four, etc. components of label vectors instead of all
components.

We start with the usual checks, then choose an empty cell's domain
and a label from it\footnote{\noindent \label{fn:branching and instantiation strategies}Heuristics
may be used to pick domains (\emph{branching strategy}) and labels
(\emph{instantiation strategy}) to optimize results. A common technique
for branching is to choose a domain among those with the least number
of labels.}, write it on the cell, then check if the resulting board is still
a PLB. If it is, we inspect the board for completion, in which case
we add it to the output list, as the board is now a solution by definition\footnote{When easy ways to derive Latin boards from a given one are available
\textendash for example isotopic or conjugacy transformations in Latin
squares, or geometric transformations in Latin polytopes\textendash{}
the algorithm could harvest instead canonical representatives of the
corresponding equivalence classes \cite{Latin Polytopes arxiv,Handbook of Constraint Programming}. }. 

If the board is indeed a PLB but still contains empty cells, we choose
a label from another cell domain, write it and perform the same checks.
We proceed in this way until the board ceases to be a PLB, at which
point we remove the just written label to restore the corresponding
empty cell, then choose a different one from its domain and resume\footnote{This is known as \emph{state restoration} \cite{Handbook of Constraint Programming}.}. 

When all labels in a cell domain have been tried, we backtrack to
the precedent cell \textendash if the current one wasn't the first\textendash{}
then, provided there are any left, choose a label that was not chosen
before and resume. If after the full traversal we have rejected all
labels in the domain we started with, we can safely conclude that
the board is not a CPLB and return the empty list.

A useful time saving device is \emph{no-good recording}\label{nogood recording}
\cite{Handbook of Constraint Programming}, a technique that remembers
combinations of values for variables found not to lead to solutions,
in order to avoid exploring again their associated search subtree\footnote{Again, when easy ways to derive Latin boards from a given one are
available, the algorithm could avoid exploring the corresponding symmetrical
no-goods \cite{Latin Polytopes arxiv,Handbook of Constraint Programming}. }. In this case backtracking may not go back to the current variable
domain but jump (\emph{back-jump}) to a different one \cite{Handbook of Constraint Programming}.

Backtracking guarantees correct termination, so when instructed to
find all solutions with\texttt{ N} $=\infty$, \texttt{list\_LBs\_from\_PB()}
will return a list with either:
\begin{enumerate}
\item no solutions, if the input PB is not a CPLB or
\item all solutions, if the input is a CPLB with several solutions or 
\item just one solution, if the input is a CPLB with just one solution\footnote{I.e.: an LP.}
\end{enumerate}

\subsubsection{Polynomial Attrition}

As the maximum count for labels $2$ and $3$ has been reached in
the bottom right sub-square of Fig. \ref{fig:sudoku ripeto 223334444}
(left), we can safely reduce the domain of each empty cell there from
$\left\{ 2,3,4\right\} $ to $\left\{ 4\right\} $. Provided that
the input PB is indeed a CPLB, $4$ is now a solution label for both
cells, and as such it can be used to trigger further reductions in
other rows and columns the cells belong to\footnote{The attrition theorem behind this reduction is immediate.}.
An attrition algorithm that visits all empty cells and applies the
theorem will be polynomial, as it will take quadratic time in the
number of cells in the board.

Additional attrition theorems with associated polynomial algorithms
may be considered\footnote{Polynomial attrition algorithms for \textit{Sudoku} may be found in
\cite{Attrition Theorems for Sudoku}.}. Properly scheduled in \texttt{list\_LBs\_from\_PB()} to maximize
effects, algorithms could parse asterisms one at a time or tuples
of them to profit from conditions arising at their intersections. 

Again, as every reduction is backed up by theorems, any cell domain
becoming empty in the process indicates that the initial assumption
for the input PB being a CPLB was not right, so \texttt{list\_LBs\_from\_PB()}
should return an empty list in this case as specified.

As before, every time a singleton domain is produced, the board should
be checked for completion. If that is the case, the solution \textendash necessarily
unique in this case\textendash{} is returned as the only element in
the output list.

A third outcome for \texttt{list\_LBs\_from\_PB()} is to stop at a
PLB that is not a solution. This will happen when the implemented
attrition algorithms cannot fully reduce the board domain.

Polynomial attrition improves then exponential attrition at the expense
of having to:
\begin{itemize}
\item prove attrition theorems
\item implement attrition algorithms
\item schedule the attrition algorithms
\item put up in some cases with non-conclusive outputs
\end{itemize}

\subsubsection{Optimized Attrition\label{sub:Optimized-Attrition}}

A way to finish the job in all cases is for \texttt{list\_LBs\_from\_PB()}
to perform exponential attrition just after polynomial attrition.
An even better option is to \emph{interleave} both techniques as we
describe next.

We start by applying polynomial attrition after the initial checks.
If this does not solve the board, we perform backtracking as described
in Sect. \ref{Backtracking}.

Any time a label is written by backtracking, and provided the resulting
board is indeed a PLB, we apply polynomial attrition again.

Every time a label is written by polynomial attrition, and in contrast
to pure polynomial attrition, we have to check now the PLB condition.
This is necessary because, once backtracking has been used, when a
label is written by the ensuing polynomial attrition the resulting
sets of label candidates for affected empty cells are no longer guaranteed
to contain all solution labels, i.e.: the sets may have ceased to
be cell domains. This has another consequence: when the board is no
longer a PLB and we have to backtrack, we have to undo all changes
made by polynomial algorithms to these \emph{pseudo-domains} since
the last label was chosen. \label{state restoration}This is known
as \emph{state restoration} \cite{Handbook of Constraint Programming}.

Along the way, we check as always for board completion, adding to
the output list every solution that may come up. As before, \texttt{list\_LBs\_from\_PB()
}is guaranteed to end with all the solutions if \texttt{N} $=\infty$,
or part thereof otherwise.

\subsection{Latin puzzles from Latin Boards\label{sub:Latin-puzzles-from-Latin Boards}}

Once \texttt{list\_LBs\_from\_PB()} has produced some LBs, we can
find LPs for each of them with Alg. \ref{pro:LPs from LB}. We start
by removing one label from the input LB, then check whether the resulting
board \textendash guaranteed now to be both a PLB and a CPLB\textendash{}
is an LP. If it is, we can try to find a more difficult puzzle by
removing another label from the just found puzzle. If the first removal
didn't produce an LP we restore the label, remove a different one,
then check again and so on. If the argument to the inner recursion
sub-algorithm in Alg. \ref{pro:LPs from LB} is passed as reference,
just a single data structure will be needed for the different board
domains and pseudo-domains generated in the process.

\subsubsection{\noindent Single Pass puzzles\label{Single-Pass-Puzzling}}

\noindent An interesting faster option to find LPs is to do so at
the same time LBs are found. This may save considerable time at the
expense of complicating \texttt{list\_LBs\_from\_PB()}. We achieve
this by both \emph{instrumenting} the attrition process and\emph{
remembering} the intermediate PLBs produced.

Given an empty PB as input, we start by finding LBs as before. Then,
any time a label is written \textendash as a consequence of either
polynomial or exponential attrition\textendash{} the operation is
uniquely tagged with a serial number or a time stamp. The tags and
their corresponding PLBs are kept in a data structure created in advance
that supports push and pop operations, that are carried out in sync
with the attrition process.

Every PLB and LB found in the process will then have a unique \emph{tag
sequence} that identifies all its antecessor PLBs. If all possibilities
for these boards have been explored, every tag in an LB's tag sequence
that does not show up in any other LB's sequence corresponds necessarily
to a puzzle. This puzzle is the stored PLB having a tagging sequence
that ends with that very tag.

\subsection{Fair Latin puzzles\label{sub:Fair-Latin-puzzles}}

A puzzle designer would like to guarantee players that her puzzles
can be solved by ``reasoning'', ``logic'', ``without trial and
error'' or ``with a pen, not a pencil''. In our terms, she would
like to produce puzzles solvable just by polynomial attrition, without
resorting to exponential methods. How can she go about it?
\begin{defn}
Given a Latin puzzle, a\emph{ polynomial attrition sequence} (PAS)
is a sequence of applications of polynomial attrition algorithms (see
Def. \ref{def:Attrition-Procedures}). A \emph{full polynomial attrition
sequence }(FPAS) is one that solves the puzzle.
\end{defn}
\noindent Let's assume our designer has solved a puzzle with a particular\emph{
}FPAS. Now, when someone else tries to solve the puzzle he may well
start using a different strategy. One involving perhaps different
polynomial algorithms, or the designer's but in a different order,
i.e: he may start using a PAS that is not a prefix of the designer's
FPAS. The question now arises: will our player be able to solve the
puzzle polynomially? Or equivalently, can any PAS be part of a FPAS
in a Latin puzzle?
\begin{defn}
A Latin puzzle is \emph{fair} if and only if any PAS is part of an
FPAS.
\end{defn}
\noindent Fair Latin puzzles are then what our designer is after.
So an equivalent question would be: under which condition is a Latin
puzzle fair? To answer, let's define first a particular class of attrition
theorems: those that achieve at least the same reduction in a domain
than in a less advanced one (see Def. \ref{def:more advanced board domains}):
\begin{defn}
Let $D$, $D'$ be two board domains for a CPLB verifying $D'\subseteq D$.
An attrition theorem $T$ is \emph{monotonic} if and only if $T\left(D'\right)\subseteq T\left(D\right)$.
\end{defn}
\noindent We also say that an attrition algorithm (see Def. \ref{def:Attrition-Procedures})
is \emph{monotonic} if it implements a monotonic attrition theorem,
and that a PAS is \emph{monotonic} if all their intervening algorithms
are monotonic.
\begin{thm}
A Latin puzzle is fair if it admits a monotonic FPAS.\end{thm}
\begin{proof}
Let $S'$ and $S$ be respectively a PAS \textendash not necessarily
monotonic\textendash{} and a monotonic FPAS for a given Latin puzzle.
Let $D\left(x_{i}\right),i=1\ldots n$, be the initial domains for
the $n$ empty cells in the puzzle. After $S'$ has been applied,
we necessarily have $D'\left(x_{i}\right)\subseteq D\left(x_{i}\right),i=1\ldots n$,
where $D'\left(x_{i}\right)$ is the new domain for cell $i$. We
apply now the first attrition algorithm in $S$. Although the cell
domains intervening in its premises may have been reduced by $S'$,
and since this algorithm is monotonic, we are certain to get a reduction
at least as strong as the one achieved if \emph{$S'$} had not been
applied. This in turns enables successive, similar application of
the rest of monotonic attrition algorithms in $S$ until the puzzle
is solved.
\end{proof}
\noindent Our designer solves then her problem by making sure that
all algorithms intervening in her own FPAS are monotonic\footnote{\noindent For more information to players, the guarantee given by
the puzzle designer could mention the set \textendash not the sequence!\textendash{}
of attrition theorems used in her FPAS.}.

\subsubsection{Finding Fair Latin puzzles}

\noindent A way of finding fair Latin puzzles is to instrument \texttt{list\_LBs\_from\_PB()}
as described in Sect. \ref{Single-Pass-Puzzling} and make it perform
optimized attrition (see Sect. \ref{sub:Optimized-Attrition}). This
time though, the data structure in the instrumentation will support
push and pop operations on attrition algorithm names instead.

Now, any time a label is written we push the name of either the polynomial
algorithm applied, or ``backtrack'' if this was the cause of the
writing. The structure is kept in sync with the attrition process
as before, so name pops will also occur when labels are removed from
the board.

At the end, every puzzle found will have an associated sequence of
logical steps towards its solution, i.e: a\emph{ proof}. The fair
puzzles among them are simply those that do not include ``backtrack''
in their proofs.

\subsubsection{Rating Fair Latin puzzles\label{sub:Rating fair Latin puzzles}}

\noindent Another advantage of puzzle proofs is that we can derive
from them objective measures of difficulty, for example as a function
of the complexity of the attrition algorithms present in the proof
and the number of times they have been applied.

This information may be used also to design a subjective scale of
difficulty, much more useful to players. To this end, we select first
a sample of players, give them samples of fair puzzles to solve and
ask for their perceived difficulty\footnote{In the scale \emph{very easy},\emph{ easy},\emph{ medium},\emph{ difficult},\emph{
very difficult} for example}. From these objective and subjective measurements, a subjective metric
for difficulty can be statistically designed. Related methods to rate
\textit{Sudoku} puzzles may be found in \cite{SudokuSat-A Tool for Analyzing Difficult Sudoku Puzzles,rating Sudoku II,rating Sudoku I}\footnote{For Latin Puzzles akin to \textit{Sudoku}, we may also collect published
\emph{Sudoku} puzzles, solve them with instrumented attrition, relate
the proofs to their featured difficulties, and use the results as
an additional ingredient in the subjective scale of difficulty.}.

\subsection{Constraint Programming}

It turns out that the problem of finding Latin boards (LBs) from a
partial labeled board (PB) is a special case of a much more general
one \cite{Principles of Constraint Programming Apt,Handbook of Constraint Programming}:
\begin{defn}
A \emph{Constraint Satisfaction Problem} (CSP) $P$ is a triple $P=\left(X,D,C\right)$
where $X$ is an $n$-tuple of variables $X=\left(x_{1},x_{2},...,x_{n}\right)$,
$D$ is a corresponding $n$-tuple of \emph{domains} $D=\left(D\left(x_{i}\right),D\left(x_{2}\right),...,D\left(x_{n}\right)\right)$,
with $x_{i}$ ranging over $D\left(x_{i}\right)$, and $C_{i}$ is
a $t$-tuple of \emph{constraints} $C=\left(C_{1},C_{2},...,C_{t}\right)$,
where a constraint $C_{j}$ is a pair $\left(R_{S_{j}},S_{j}\right)$
in which $R_{S_{j}}$ is a relation on the variables in $S_{i}=$
\emph{scope}$(C_{i})$.
\end{defn}
\noindent A constraint is then a subset of the Cartesian product of
its variables. If we identify the tuple of empty cells in a partial
labeled board with $X$, a board domain with $D$, the cell domains
with $D\left(x_{i}\right)$, and the need for every asterism in the
board to hold all labels in the multiset with a constraint, then the
problem to find solution Latin boards for a partial labeled board
can be formulated as a Constraint Satisfaction Problem (CSP).

\subsubsection{Constraint Programming Systems (CPS)}

This new way of putting things does spare many efforts when specifying
and implementing algorithm \texttt{list\_LBs\_from\_PB()} (see Alg.
\ref{pro:LBs from PB}), as it also turns out that many of the services
it should provide are readily available in \emph{Constraint Programming
Systems}\footnote{Or more specifically, as our variable domains are finite\emph{, }in\emph{
Finite Domain CPSs} \cite{Handbook of Constraint Programming}.} (CPS) \cite{Principles of Constraint Programming Apt,Handbook of Constraint Programming}.

The main feature of a CPS is that the user simply \emph{declares}
the constraints in some language \textendash usually a subset of first-order
logic\textendash{} then lets the system run to find the answers. This
programming paradigm differs notably from imperative ones, as the
focus here is not the solving algorithm, but rather the specification
of the properties of the solution.

A CPS can be used to tackle a wide range of combinatorial optimization
problems, for example that of completing partial Latin squares \cite{Completing quasigroups or latin squares,Arc consistency and quasigroup completion}.
A suitable CPS can also be used to solve \textit{Sudoku} puzzles,
both alone \cite{Using Constraint Programming to solve Sudoku Puzzles,Overlapping Alldifferent constraints and the Sudoku puzzle,Generating and Solving Logic Puzzles through Constraint Satisfaction,Sudoku as a constraint problem}
or in conjunction with other techniques \cite{Nonrepetitive paths and cycles in graphs with application to Sudoku,Combining Metaheuristics and CSP Algorithms to Solve Sudoku,Building Industrial Applications with Constraint Programming,A Hybrid alldifferent-Tabu Search Algorithm for Solving Sudoku Puzzles}.

In a CPS, attrition theorems (see Def. \ref{sub:Attrition-Theorems})
are called \emph{rules }or\emph{ proof rules. }Implemented polynomial
attrition algorithms (see Def. \ref{def:Attrition-Procedures}) are
called \emph{propagators};\emph{ }their application, \emph{constraint
propagation} or \emph{inference}, and exponential attrition (see Sect.
\ref{sub:Types-of-Attrition}),\emph{ search. }When scheduled by a
\emph{propagation engine}\label{Propagation engine}, a CPS performs\emph{
}propagation\emph{ }by maintaining\emph{ network consistency}, which
means carrying out safe attrition in all its forms. CPSs come with
propagators that enforce \emph{local consistency} for individual variables
(\emph{node consistency}), pairs of variables (\emph{arc consistency})
and tuples of variables (\emph{path consistency}) with the aim of
obtaining \emph{global network consistency} \cite{Handbook of Constraint Programming}.

CPSs are implemented either as autonomous systems or libraries \cite{Handbook of Constraint Programming},
and frequently allow users to plug in their own propagators and search
strategies. In fact, the author used a custom CPS to generate the
example Latin boards and puzzles featured in Sect. \ref{sec:Examples of Latin puzzles}.

\subsection{Alternatives to Constraint Programming \label{sub:Alternatives-to-CP}}

A CPS equipped with suitable propagators is a good tool to find and
solve Latin puzzles, but there are alternatives. We could in principle
consider other combinatorial optimization techniques known to solve
partial Latin squares and \textit{Sudoku} puzzles, and adapt them
to work with the multisets and arbitrary topologies proper to Latin
puzzles. Some of these techniques are\footnote{For related techniques specific to combinatorial designs see Sect.
\emph{Computational Methods in Design Theory} in \cite{Handbook of Combinatorial Designs}.}:
\begin{itemize}
\item propositional satisfiability inference (SAT) \cite{Sudoku as a SAT problem}
\item Gröbner bases \cite{Gr=0000F6bner Basis Representations of Sudoku,Gr=0000F6bner bases,Sudokus and Gr=0000F6bner Bases}
\item rewriting logic \cite{Solving Sudoku Puzzles with Rewriting Rules}
\item dancing links \cite{Dancing links}
\item another solution problem \cite{Complexity and completeness of finding another solution and its application to puzzles}
\item graphs and hypergraphs coloring \cite{Coloring Mixed Hypergraphs}
\item ant colony optimization \cite{Solving sudoku puzzles using hybrid ant colony optimization algorithm,An Ant Algorithm for the Sudoku Problem}
\item simulated annealing \cite{Metaheuristics can solve sudoku puzzles}
\item tabu search \cite{Tabu search}
\item genetic algorithms \cite{Solving Sudoku Puzzles with Wisdom of Artificial Crowds,Solving and analyzing Sudokus with cultural algorithms,Product geometric crossover for the sudoku puzzle,Solving sudoku with the gAuGE system}
\item bee colony simulation \cite{Solving sudoku puzzles using improved artificial bee colony algorithm}
\item harmony search \cite{Harmony search algorithm for solving sudoku}
\item flower pollination \cite{A Novel Hybrid Flower Pollination Algorithm with Chaotic Harmony Search for Solving Sudoku Puzzles}
\item particle swarm optimization \cite{Geometric particle swarm optimization for the sudoku puzzle}
\item entropy measurement \cite{Entropy Minimization for Solving Sudoku,Solving Sudoku Puzzles Based on Customized Information Entropy}
\end{itemize}
Deterministic techniques like SAT, Gröbner bases or dancing links,
which are able in principle to find all solutions to a partial labeled
board\footnote{We saw in Sect. \ref{sub:Latin-puzzles-from-Latin Boards} that this
feature is critical to find Latin puzzles.}, could be used to find Latin puzzles, while stochastic methods would
be suitable to solve them.

\section{Conclusions}

The generalization of Latin squares to Latin boards proposed in \cite{Latin Polytopes arxiv}
prompts new lines of research. In this paper we have suggested a corresponding
generalization of Latin square puzzles, together with methods to obtain
Latin puzzles and examples thereof. Future work could proceed along
the following lines:
\begin{itemize}
\item new Latin boards and Puzzles, in addition to the ones featured in
this paper and those in \cite{Latin Polytopes arxiv,Latin Puzzles site}
\item monotonic attrition theorems and propagators to solve partial Latin
boards (see Sect. \ref{sub:Latin-Boards-from-Partial-Boards})
\item techniques to solve partial Latin boards other than constraint programming
(see Sect. \ref{sub:Alternatives-to-CP})
\item applications in scheduling, timetabling and resource allocation problems
when solutions are Latin boards
\item generalization of existing results for partial Latin squares and quasi-groups
in Algebra:

\begin{itemize}
\item notions of isotopism and conjugacy in partial Latin boards (for partial
Latin squares see \cite{The set of autotopisms of partial Latin squares})
\item connection with other branches of discrete mathematics \cite{Sudoku gerechte designs resolutions affine space spreads reguli and Hamming codes}
\item enumeration of particular Latin boards (for \textit{Sudoku} see \cite{Enumerating possible Sudoku grids})
\item number of clues in minimal Latin puzzles, with and without inscription
(for \textit{Sudoku} see \cite{There is no 16-Clue Sudoku,Sudoku dataset: 49151 instances with 17 givens})
\item NP-completeness of the problem of completing specific partial Latin
boards, like partial Latin triangles, partial Latin hexagons, etc.
(for partial Latin squares see \cite{Completing partial LS is NP-complete})
\end{itemize}
\end{itemize}

\section*{Acknowledgements\addcontentsline{toc}{section}{Acknowledgements}}

The author thanks Latin squares expert Prof. Raúl Falcón Ganfornina
and puzzle expert Will Shortz for their reviews of this paper, as
well as the following persons for their comments: Francisco Javier
Acinas Lope, Roberto Álvarez Chust, Prof. Francisco Ballesteros Olmo,
Ceit Hustedde, Marie-Aimée Lehen, Jesús León Álvaro Páez, Eduardo
Pérez Molina, Fe Pesquero Martínez, Andrew Thean and Manuel Valderrama
Conde. The author also thanks the members of \textit{Tertulia de Matemáticas}
\cite{Tertulia de Matem=0000E1ticas} for all the great conversations.

\pagebreak{}

\tableofcontents{}\addcontentsline{toc}{section}{Contents}

\pagebreak{}

\listoffigures

\addcontentsline{toc}{section}{List of Figures}

\bigskip{}

\paragraph{About the author}

\noindent \textbf{\addcontentsline{toc}{section}{About the author}}{\small{}The
author is a telecommunication engineer from the Polytechnical University
of Madrid. He has created several mathematical Puzzles like }\emph{\small{}Sudoku
Ripeto}{\small{} \cite{Latin Puzzles site}, }\emph{\small{}Custom
Sudoku }{\small{}\cite{Latin Puzzles site}}\emph{\small{},}{\small{}
}\emph{\small{}Konseku}{\small{} \cite{Konseku Puzzle site} and}\emph{\small{}
Moshaiku}{\small{} \cite{Moshaiku Puzzle site}.}{\small \par}
\end{document}